# The Complexity of Approximating PSPACE-Complete Problems for Hierarchical Specifications [1,2,3]

Madhav V. Marathe    Harry B. Hunt III    S.S. Ravi

Department of Computer Science

University at Albany - State University of New York

Albany, NY 12222

## Abstract

We extend the concept of polynomial time approximation algorithms to apply to problems for hierarchically specified graphs, many of which are PSPACE-complete. Assuming P $\neq$ PSPACE, the existence or nonexistence of such efficient approximation algorithms is characterized, for several standard graph theoretic and combinatorial problems. We present polynomial time approximation algorithms for several standard PSPACE-hard problems considered in the literature. In contrast, we show that unless P = PSPACE, there is no polynomial time $\epsilon$-approximation for any $\epsilon > 0$, for several other problems, when the instances are specified hierarchically.

We present polynomial time approximation algorithms for the following problems when the graphs are specified hierarchically:

*minimum vertex cover, maximum 3SAT, weighted max cut, minimum maximal matching,* and

*bounded degree maximum independent set.*

In contrast, we show that unless P = PSPACE, there is no polynomial time $\epsilon$-approximation for any $\epsilon > 0$, for the following problems when the instances are specified hierarchically:

*the number of true gates in a monotone acyclic circuit when all input values are specified* and

*the optimal value of the objective function of a linear program.*

It is also shown that unless P = PSPACE, a performance guarantee of less than 2 cannot be obtained in polynomial time for the following problems when the instances are specified hierarchically:

*high degree subgraph, k-vertex connected subgraph* and *k-edge connected subgraph.*

**Classification**: Hierarchical Specifications, Approximation Algorithms, Computational Complexity, Algorithms and Data structures.

[1] Email addresses of authors:{madhav,hunt,ravi}@cs.albany.edu

[2] Supported by NSF Grants CCR 89-03319 and CCR 89-05296.

[3] An extended abstract of this paper appeared in the proceedings of *20th International Colloquium on Automata, Languages and Programming (ICALP)*, 1993.

# 1 Introduction

Hierarchical system design is becoming increasingly important with the development of VLSI technology [HLW92, RH93]. At present, a number of VLSI circuits already have over a million transistors. (For example the *Intel i860* chip has about 2.5 million transistors.) Although VLSI circuits can have millions of transistors, they usually have highly regular structures. These regular structures often make them amenable to hierarchical design, specification and analysis. Other applications of hierarchical specification and consequently of hierarchically specified graphs are in the areas of finite element analysis [LW87b], software engineering [GJM91], material requirement planning and manufacturing resource planning in a multistage production system [MTM92] and processing hierarchical Datalog queries [Ul88].

Over the last decade, several theoretical models have been put forward to succinctly represent objects hierarchically [BOW83, GW83, Le82, Le88, LW93, Wa86]. Here, we use the model defined by Lengauer in [HLW92, Le86, Le88, LW92] to describe graphs. Using this model, Lengauer et al. [Le89, LW87a, Le88] have given efficient algorithms to solve several graph theoretic problems including minimum spanning forests, planarity testing etc.

Here, we extend the concept of polynomial time approximation algorithms so as to apply to problems for hierarchically specified graphs including PSPACE-complete such problems. We characterize the existence or nonexistence (assuming P $\neq$ PSPACE) of polynomial time approximation algorithms, for several standard graph problems. Both positive and negative results are obtained (see Tables 1 and 2 at the end of this section). Our study of approximation algorithms for hierarchically specified problems is motivated by the following two facts:

1. $\Theta(n)$ size hierarchical specifications can specify $2^{\Omega(n)}$ size graphs.

2. Many basic graph theoretic properties are PSPACE-complete [HR+93, LW92], rather than NP-complete.

For these reasons, the known approximation algorithms in the literature are not directly applicable to graph problems, when graphs are specified hierarchically.

What we mean by a *polynomial time approximation algorithm* for a graph problem, when the graph is specified hierarchically, can be best understood by means of an example.

**Example:** Consider the minimum vertex cover problem, where the input is a hierarchical specification of a graph $G$. We provide efficient algorithms for the following versions of the problem.

1. **The Approximation Problem:** Compute the size of a near-minimum vertex cover of $G$.

2. **The Query problem**: Given any vertex $v$ of $G$ and the path from the root to the node in the *hierarchy tree* (see Section 2 for the definition of hierarchy tree) in which $v$ occurs, determine whether $v$ belongs to the approximate vertex cover so computed.

3. **The Construction Problem:** Output a hierarchical specification of the set of vertices in the approximate vertex cover.

4. **The Output Problem:** Output the approximate vertex cover computed.

Our algorithms for (1), (2) and (3) above run in time **polynomial in the size of the hierarchical specification** rather than the size of the graph obtained by expanding the specification. Our algorithm for (4) runs in time linear in the size of the expanded graph but uses space which is linear in the size of the hierarchical specification. ∎

This is a natural extension of the definition of approximation algorithms in the flat (i.e. non-hierarchical) case. This can be seen as follows. In the flat case, the number of vertices is polynomial in the size of the description. Given this, any polynomial time algorithm to determine if a vertex $v$ of $G$ is in the approximate minimum vertex cover can be modified easily into a polynomial time algorithm that lists all the vertices of $G$ in the approximate minimum vertex cover. For an optimization problem or a query problem, our algorithms use space and time which are low level polynomials in the size of the hierarchical specification and thus $O(poly \log \eta)$ in the size of the specified graph, when the size $\eta$ of the graph is exponential in the size of the specification. Moreover, when we need to output the subset of vertices, subset of edges, etc. corresponding to a vertex cover, maximal matching, etc., in the expanded graph, our algorithms take essentially the same time but substantially less (often exponentially less) space than algorithms that work directly on the expanded graph. It is important to design algorithms which work directly on the hierarchical specification by exploiting the regular structure of the underlying graphs, because, graphs resulting from expansions of given hierarchical descriptions are frequently too large to fit into the main memory of a computer [Le86]. This results in a large number of page faults while executing the known algorithms on the expanded graph. Hence, standard algorithms designed for flat graphs are impractical for hierarchically specified graphs.

We believe that this is the first time efficient approximation algorithms with good performance guarantees have been provided both for hierarchically specified problems and for PSPACE-complete problems.[4] Thus by providing algorithms which exploit the underlying structure, we extend the range of applicability of standard algorithms so as to apply to a much larger set of instances. Tables 1 and 2 summarize our results.

---

[4] Independently, Condon et al. [CF+93a, CF+93a] have investigated the approximability of other PSPACE-complete problems.

**Table 1. Performance Guarantees**

| Problem | Performance guarantee in hierarchical case | Best known guarantee in flat case |
|---|---|---|
| MAX 3SAT | 2 | 4/3 |
| MIN Vertex Cover | 2 | 2 |
| MIN Maximal Matching | 2 | 2 |
| Bounded Degree (B) MAX Independent Set | $B$ | $B$ |
| MAX CUT | 2 | 2 |

The results mentioned in the last column of the above table can be found in [GJ79, Ya92].

**Table 2. Hardness Results**

| Problem | Hierarchical Case | Flat Case |
|---|---|---|
| Maximum Number of True Gates in a circuit | PSPACE-hard to approximate for any $\epsilon$ | Log-hard for P to approximate for any $\epsilon$ |
| Optimal Value of Objective Function of a Linear Program | PSPACE-hard for any $\epsilon$ | Log-hard for P to approximate for any $\epsilon$ |
| High Degree Subgraph | PSPACE-hard for $\epsilon < 2$ | Log-hard for P to approximate for $\epsilon < 2$ |
| $k-$ Vertex Connectivity | PSPACE-hard for $\epsilon < 2$ | Log-hard for P to approximate for $\epsilon < 2$ |
| $k-$ Edge Connectivity | PSPACE-hard for $\epsilon < 2$ | Log-hard for P to approximate for $\epsilon < 2$ |

The results mentioned in the last column of the above table can be found in [AM86, KSS89, Se91].

## 2 Definitions and Description of the Model

The following two definitions are from Lengauer [Le89].

**Definition 2.1** *A hierarchical specification $\Gamma = (G_1, ..., G_n)$ of a graph is a sequence of undirected simple graphs $G_i$ called cells. The graph $G_i$ has $m_i$ edges and $n_i$ vertices. $p_i$ of the vertices are distinguished and are called pins. The other $(n_i - p_i)$ vertices are called inner vertices. $r_i$ of the inner vertices are distinguished and are called nonterminals. The $(n_i - r_i)$ vertices are called terminals.*

Note that there are $n_i - p_i - r_i$ vertices defined explicitly in $G_i$. We call these *explicit vertices*. Each pin of $G_i$ has a unique label, its *name*. The pins are assumed to be numbered from 1 to $p_i$. Each nonterminal in $G_i$ has two labels, a *name* and a *type*. The type is a symbol from $G_1, ..., G_{i-1}$. If a nonterminal vertex $v$ is of the type $G_j$, then the terminal vertices which are the neighbors of $G_j$ are in one-to-one correspondence with the pins of $G_j$. (Note that all the neighbors of a nonterminal vertex must be terminals. Also, a terminal vertex may be a neighbor of several nonterminal vertices.) The size of $\Gamma$, denoted by $size(\Gamma)$, is $N + M$, where the vertex number $N = \sum_{1 \le i \le n} n_i$, and the edge number $M = \sum_{1 \le i \le n} m_i$.

**Definition 2.2** *Let $\Gamma = (G_1, ..., G_n)$ be a hierarchical specification of a graph $G$. The expansion $E(\Gamma)$ (i.e. the graph associated with $\Gamma$) of the hierarchical specification $\Gamma$ is done as follows:*
$k = 1 : E(\Gamma) = G_1$.
$k > 1$ : *Repeat the following step for each nonterminal $v$ of $G_k$, say of the type $G_j$: delete $v$ and the edges incident on $v$. Insert a copy of $E(\Gamma_j)$ by identifying the $l^{th}$ pin of $E(\Gamma_j)$ with the node in $G_k$ that is labeled $(v, l)$. The inserted copy of $E(\Gamma_j)$ is called the subcell of $G_k$. (Observe that the expanded graph can have multiple edges although none of the $G_i$ have multiple edges.)*

The expansion $E(\Gamma)$ is the graph associated with the hierarchical definition $\Gamma$. Note that the total number of nodes in $E(\Gamma)$ can be $2^{\Omega(N)}$. For $1 \le i \le n$, $\Gamma_i = (G_1, ..., G_i)$ is the hierarchical specification of the graph $E(\Gamma_i)$. Given a hierarchical specification $\Gamma$, one can associate a natural tree structure depicting the sequence of calls made by the successive levels. We call it the *hierarchy tree* and denote it by $HT(\Gamma)$. A vertex in $E(\Gamma)$ is identified by a sequence of nonterminals on the path from the root to the nonterminal in which the vertex is explicitly defined. For the query problems considered in the paper, we assume that a vertex is specified in the above manner.

Without loss of generality we assume that there are no useless cells in $\Gamma_n = \Gamma$.

**Example:** Figure 1 shows an example of a hierarchically specified graph and its corresponding hierarchy tree. The labels on the vertices are omitted and the 1-1 correspondence between the pins of $G_j$ and the neighbors of a nonterminal of type $G_j$ in the cell $G_i$ is clear by the positions of the vertices in the figure. Figure 2 shows the underlying graph $E(G)$. *We note again that our approximation algorithms answer query problems without explicitly expanding the hierarchical specification.*

**Definition 2.3** *A hierarchical graph specification $\Gamma = (G_1, ..., G_n)$ of a graph $G$ is **1-level-restricted** if for all $(u, v) \in E$, one of the following conditions holds :*

1. *$u$ and $v$ are explicit vertices in the same instance of $G_i$ ($1 \le i \le n$).*

2. *$u$ is an explicit vertex in an instance of $G_i$ and $v$ is a explicit vertex in an instance of $G_j$ and the instance of $G_i$ calls the instance of $G_j$ ($1 \le j < i \le n$).*

A hierarchical graph specification $\Gamma = (G_1, ..., G_n)$ of a graph $G$ is **strongly 1-level-restricted** if it is 1-level-restricted and in addition for $2 \le i \le n$, the only nonterminals of $G_i$ are of the type $G_{i-1}$.

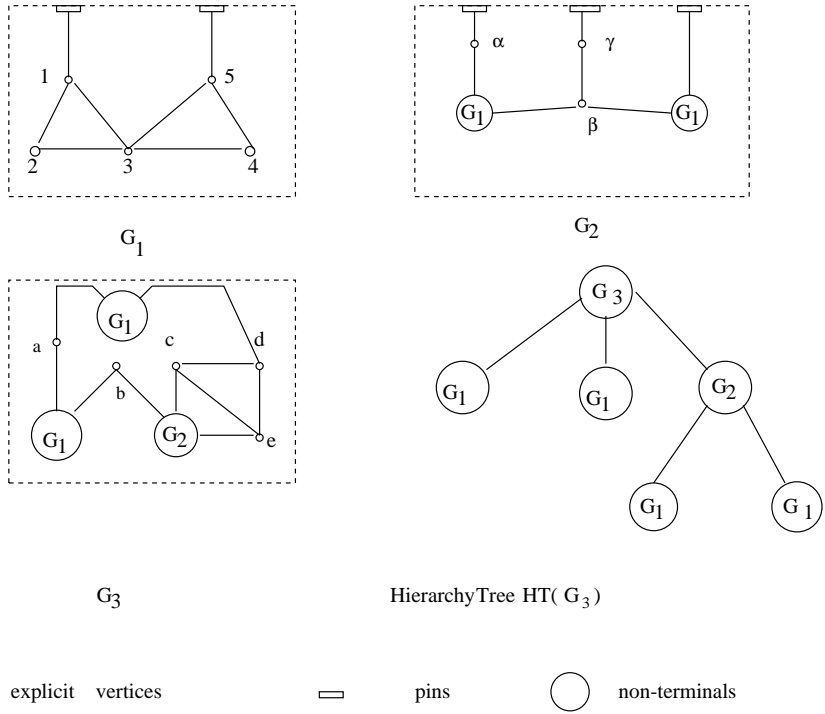

Figure 1: A hierarchically specified graph, and the corresponding hierarchy tree.

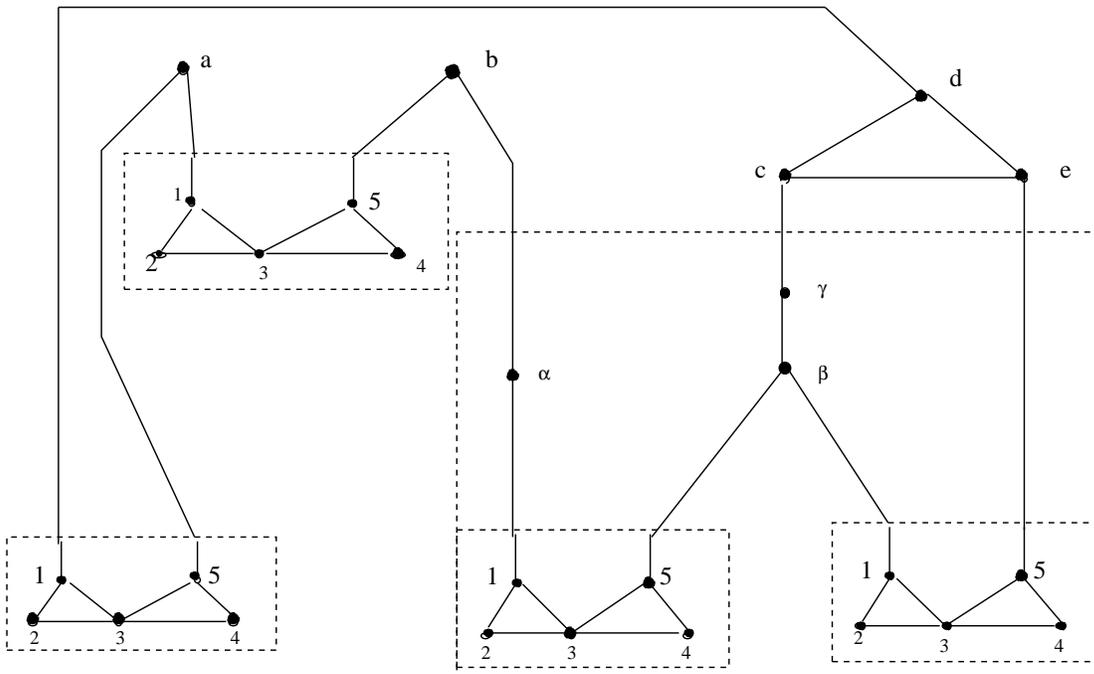

Figure 2: The graph associated with the hierarchical specification $G$.

The above definition can be extended to define *k-level restricted* specifications, for any fixed $k \geq 1$. Such descriptions still can lead to exponentially large graphs. Moreover, many practically occurring hierarchical descriptions (see [Le82, Le86, LW87a]) are *k-level restricted* for small values of $k$. We note that our PSPACE-hardness results hold for strongly 1-level-restricted specifications, while all our approximation algorithms hold for arbitrary specifications.

**Definition 2.4** *Let $\Gamma = (G_1, ..., G_n)$ be a hierarchical specification. $\Gamma$ is said to be **simple** if, for each $G_i$, $1 \leq i \leq n$, there are no edges between pins defined in $G_i$.*

For a simple 1-level-restricted specifications observe that:

**Observation 2.5** *Consider any edge $(u, v)$ in a simple 1-level-restricted hierarchical specification of a graph $G$. Then the path from $u$ to $v$ in the hierarchy tree passes through at most one pin.* ∎

For the rest of the discussion, given a problem $\Pi$ we denote by $\Pi_{HG}$ the same problem when the instance is specified hierarchically. So for example, we use $MAXCUT_{HG}$ to denote the MAX CUT problem when the graph is specified hierarchically. Also, we sometimes use the phrase hierarchical graphs to mean hierarchically specified graphs.

Finally, we give additional definitions used in the paper.

**Definition 2.6** *The **Monotone Circuit Value Problem** $MVCP$ is defined as follows: Given an acyclic graph $G$ called the circuit with one distinguished vertex (output), the sources (inputs) labeled with $\{0,1\}$ and all other vertices labeled with symbols from $\{\vee, \wedge\}$, the decision version of the problem asks if the output of $G$ is 1. The optimization version of $MCVP$ denoted by $MTG$ asks for the maximum number of gates which are set to 1.*

We assume that the reader is familiar with the problem 3SAT. The problem $3SAT_{HG}$ is defined as follows:

**Definition 2.7** *An instance $F = (F_1(X^1), \ldots, F_{n-1}(X^{n-1}), F_n(X^n))$ of $3SAT_{HG}$ is of the form*

$$F_i(X^i) = (\bigwedge_{1 \leq j \leq l_i} F_{i_j}(X^i_j, Z^i_j)) \bigwedge f_i(X^i, Z^i)$$

*for $1 \leq i \leq n$ where $f_i$ are 3CNF formulae, $X^n = \phi$, $X^i, X^i_j, Z^i, Z^i_j$, $1 \leq i \leq n-1$, are vectors of boolean variables such that $X^i_j \subseteq X^i$, $Z^i_j \subseteq Z^i$, $0 \leq i_j < i$. Thus, $F_1$ is just a 3CNF formula. An instance of $3SAT_{HG}$ specifies a 3CNF formula $f$, that is obtained by expanding the $F_j$, $2 \leq j \leq n$ as macros where the variables $Z$'s introduced in any expansion are considered distinct. The problem $3SAT_{HG}$ is to decide whether the formula $f$ specified by $F$ is satisfiable. The optimization problem denoted by MAX $3SAT_{HG}$ is to find an assignment to the variables of $f$ satisfying the maximum number of clauses in $f$.*

Let $n_i$ be the total number of variables used in $F_i$ (i.e. $|X^i| + |Z^i|$) and let $m_i$ be the total number of clauses in $F_i$. The size of $F$, denoted by $size(F)$, is equal to $\sum_{1 \leq i \leq n}(m_i + n_i)$.

**Example:** Let $F = (F_1(x_1, x_2), F_2(x_3, x_4), F_3)$ be an instance of 3SAT$_{HG}$ where each $F_i$ is defined as follows:

$$F_1(x_1, x_2) = (x_1 + x_2 + z_1) \wedge (z_2 + z_3)$$

$$F_2(x_3, x_4) = F_1(x_3, z_4) \wedge F_1(z_4, z_5) \wedge (z_4 + z_5 + x_4)$$

$$F_3 = F_2(z_8, z_7) \wedge F_1(z_7, z_6)$$

The formula $f$ denoted by $F$ is $(z_7 + z_6 + z_1^1) \wedge (z_2^1 + z_3^1) \wedge (z_8 + z_4 + z_1^2) \wedge (z_2^2 + z_3^2) \wedge (z_4 + z_5 + z_1^3) \wedge (z_2^3 + z_3^3) \wedge (z_4 + z_5 + z_7)$.

**Definition 2.8** *Let $F$ be an instance of the problem 3SAT with set of variables $V$ and set of clauses $C$.*

*1. The bipartite graph of $F$, denoted $\mathbf{BG}(f)$, is the bipartite graph $(V \cup C, E)$, where $e = (c, v) \in E$ iff variable $v$ occurs in clause $c$.*

*2. $F$ is said to be planar iff the graph $\mathbf{BG}(f)$ is planar.*

**Definition 2.9** *An instance $F = (F_1(X^1), \ldots, F_{n-1}(X^{n-1}), F_n(X^n))$ $(\Delta_1, \ldots, \Delta_{n-1}, \Delta_n)$ of Hierarchical Linear Program (LP$_{HG}$) is of the form*

$$F_i(X^i) = (\bigcup_{1 \leq i_j \leq i} F_{i_j}(X_j^i, Z_j^i)) \bigcup f_i(X^i, Z^i)$$

$$\Delta_i = \sum_{i_j} d_{i_j} \cdot \Delta_{i_j} + \sum_{z_j \in Z^i} c_j \cdot z_j$$

*for $1 \leq i \leq n$ where $f_i$ is a set of linear inequalities, $X^n = \phi$, $X_j$, $X_j^i$, $Z^i, Z_j^i$, $1 \leq i \leq n-1$, are vectors of variables such that $X_j^i \subseteq X^i$, $Z_j^i \subseteq Z^i$, $1 \leq i_j \leq i$, $F_i$ is a set of linear inequalities and $\Delta_i$ is a linear objective function over the variables in $E(F_i)$. Thus $F_1$ is just a set of linear inequalities. An instance of LP$_{HG}$ defines a hierarchically specified linear program $F_n$ obtained after expanding $F_j$ ($1 \leq j \leq n$) as macros where the Z's in different expansions are considered distinct and a linear objective function $\Delta_n$ obtained after expanding $\Delta_j's$ as macros.*

Let $n_i$ be the total number of variables used in $F_i \cup \Delta_i$ and let $m_i$ be the total number of inequalities in $F_i$. Then, the size of $F$ denoted by $size(F)$ is equal to $\sum_{1 \leq i \leq n}(m_i + n_i)$.

The $LP$ feasibility problem is to determine whether there exists an assignment to the variables (over the reals) used in the $LP$, such that all the inequalities are satisfied. In the case of the LP$_{HG}$ optimization problem, one is given a linear objective function and linear inequalities both defined hierarchically as above. The aim is to find an assignment to the variables so as to maximize the value of the objective function subject to the inequality constraints. Using Lengauer's definition of hierarchical graphs, one can represent a $LP_{HG}$ graphically by associating a node with each variable and with each inequality. Further, a variable node has an edge to an inequality node iff the corresponding variable occurs in the inequality.

Linear programming has been extensively studied in literature. In [GLS84] it is shown how linear programs can be used to model many graph theoretic problems. In [GLS84] it was also shown that for the class of *perfect graphs*, polynomial time algorithms can be devised to compute an optimal vertex coloring, maximum independent set and several other important graph theoretic parameters. When graphs are represented hierarchically, the corresponding linear program will be hierarchical. But as will be shown (Section 7), computing the optimal value of the objective function of a hierarchically specified linear program is PSPACE-hard; further, it is also PSPACE-hard to compute an approximate value of the objective function.

Next, we recall the definitions of high degree subgraph and high vertex (edge) connectivity problems.

**Definition 2.10** *The* **High Degree Subgraph Problem ($k$-HDSP)** *is defined as follows: For all integers $k \geq 3$, given a graph $G = (V, E)$, does $G$ have a nonempty subgraph of minimum degree $k$. The optimization problem of $k$-HDSP, denoted by MAX HDSP, asks for the maximum $k$ such that there is a vertex induced subgraph of $G$ in which the minimum degree of a vertex is $k$.*

Let HDSP$^*$ denote the largest $k$ such that there is an induced subgraph of minimum degree $k$. An approximate solution to this problem is a subgraph in which each node has degree at least $d$, where HDSP$^* \geq d \geq$ HDSP$^*/c$, for some fixed $c > 1$. For all $k \geq 3$, $k$-HDSP was shown to be log-complete for $P$ in [AM86]. Furthermore, unless $P = NC$, it was shown that no NC approximation algorithm for MAX HDSP could provide a performance guarantee better than 2. $k$-HDSP is polynomial time solvable for flat graphs [AM86]. We show that $k$-HDSP$_{HG}$ is PSPACE-complete and furthermore unless P = PSPACE, MAX HDSP$_{HG}$ cannot be approximated with a factor $c < 2$ in polynomial time (See Section 7). The high degree subgraph problem contrasts with the related *maximum clique problem (MCP)* which is NP-complete for both flat [GJ79] and hierarchically specified graphs [LW92].

Next we recall the definitions the high-vertex and edge connectivity problems from [KSS89].

**Definition 2.11** *The vertex connectivity $\kappa(G)$ (edge connectivity $\lambda(G)$) of an undirected graph $G$ is the minimum number of vertices (edges) whose removal results in a disconnected or a trivial graph[5] A graph is m-vertex-connected (m-edge-connected) if $\kappa(G) \geq m$ ($\lambda(G) \geq m$).*

**Definition 2.12** *The* **High Vertex Connectivity Problem ($\kappa$-HVCP)** *(*High Edge Connectivity Problem ($\kappa$-HECP)*) is defined as follows: For all integers $\kappa \geq 3$, given a graph $G = (V, E)$, does $G$ contain an induced subgraph of vertex connectivity (edge connectivity) at least $\kappa$? The optimization*

---
[5]A trivial graph consists solely of isolated vertices.

*versions of these problems denoted by MAX HVCP (MAX HECP) ask for the largest $\kappa$ for such that there is an induced subgraph of vertex(edge) connectivity $\kappa$.*

Let HVCP* (HECP*) denote the largest $\kappa$ such that there is an induced subgraph of vertex(edge) connectivity $\kappa$. An approximate solution to this problem is a subgraph whose vertex (edge) connectivity is at least $d$, where HVCP* (HECP*) $\geq d \geq$ HVCP*/$c$ (HECP*/$c$), for some fixed $c > 1$. It was shown in [KSS89] that for all $\kappa \geq 3$, $\kappa$-HVCP and $\kappa$-HECP are log-complete for $P$. Furthermore, they showed that

**Theorem 2.13** (Kirousis, Serna, Spirakis [KSS89]) *Unless $P \neq NC$, MAX HVCP and MAX HECP cannot be approximated to within a factor $c < 2$ of the optimal in NC.*

Here, we show that for all $\kappa \geq 3$, the problems $\kappa$-HVCP$_{HG}$ and $\kappa$-HECP$_{HG}$ are PSPACE-complete and furthermore unless P = PSPACE, MAX HVCP$_{HG}$ and MAX HECP$_{HG}$ cannot be approximated within a factor $c < 2$ in polynomial time (See Section 7).

We end this section with a few comments regarding our approximation algorithms for the problems MAX-CUT, MAX 3SAT and Bounded-degree Independent set when instances are specified hierarchically. Consider the MAX CUT problem. For any graph $G(V, E)$, there is always a cut containing at least $|E|/2$ edges. Therefore, by merely counting the number of edges in a hierarchically specified graph, one can always compute a number which is within a factor of 2 of an optimal cut. However our approximation algorithm for the MAX CUT problem *actually finds* a hierarchical representation of a cut containing at least $|E|/2$ edges. Similar comments apply to our approximation algorithms for the problems MAX 3SAT and Bounded-degree Independent set when instances are specified hierarchically.

# 3 Approximation Algorithms

In this section we discuss our approximation algorithms for the problems given in Table 1. We first outline the basic technique used to efficiently obtain approximation algorithms with good performance guarantee.

## 3.1 The Basic Technique: Approximate Burning

Our approximation algorithms are based on a new technique which we call **approximate burning**. This is an extension of the **Bottom Up** method for processing hierarchical graphs discussed in [LW87a, Le88, Le89] and [Wi90] for designing efficient algorithms for hierarchically specified graphs. The bottom up method aims at finding a small graph $G_i^b$ called the *burnt graph* which can replace each occurrence of $G_i$ in such a way that $G_i$ and $G_i^b$ behave identically with respect to the problem under consideration. The bottom up method should produce such burnt graphs efficiently. Since the problems we are dealing with are PSPACE-hard, we cannot hope to find in polynomial time such burnt graphs which can replace original graphs. Therefore, we resort to *approximate burning*. In approximate burning, given an approximation algorithm for non-hierarchical instances of the problem, we wish to find small burnt graphs which can be used to replace the original non-terminals in such a way that the performance guarantee provided by the algorithm is not affected by the replacement. All our approximation algorithms rely on approximate burning.

In summary, to obtain good solutions for a problem specified hierarchically, the bottom up procedure should have the following properties:

1. Each burnt graph should have a size which is polynomial in the size of the specification.

2. The burning procedure should run in time which is polynomial in the size of specification.

3. The burnt graphs should be replaceable with respect to the problem $\Pi$ and the approximation algorithm $A_\Pi$.

Before we discuss our approximation algorithm, we give a transformation which allows us to transform a hierarchical specification in which there are edges between pins defined in $G_i$ to an equivalent hierarchical specification which has no edges between pins defined in a given $G_i$. The transformation is outlined in Figure 3.

The following lemma summarizes the property of the specification $\Gamma_1$ obtained as a result of the transformation outlined in Figure 3.

**Lemma 3.1** *Given a hierarchical specification $\Gamma = (G_1, ..., G_n)$ in which there are edges between pins defined in a given $G_i$, we can construct in polynomial time n new hierarchical specification $\Gamma_1 = (H_1, ..., H_n)$ such that*

1. *$size(\Gamma_1)$ is polynomial in $size(\Gamma)$.*

2. *$\Gamma_1$ can be constructed in polynomial time.*

3. *$E(\Gamma) = E(\Gamma_1)$.*

4. *For each $H_i$, $1 \leq i \leq n$, there are no edges between pins defined in $H_i$.* ∎

In view of Lemma 3.1, we assume that in the input to all our approximation algorithms is a *simple* hierarchical specification (i.e. there is no edge between two pins which are defined in the same cell). The running times of our approximation algorithms are with respect to such simple specification.

### 3.2 Approximation Algorithm for Vertex Cover

We now discuss our heuristic for computing the size of a near-optimal vertex cover for a hierarchically specified graph. The problem of computing the size of a minimum vertex cover for hierarchically specified graphs was shown to be PSPACE-hard by Lengauer [LW92] (Actually, they prove the hardness for maximum independent set; the hardness of minimum vertex cover is therefore directly implied). Our heuristic

**Procedure Transform-HSPEC**
**Input:** *A hierarchical specification* $\Gamma = (G_1, ..., G_n)$ of a graph $G$.
**Output:** *A new hierarchical specification* $\Gamma_1 = (H_1, ..., H_n)$ which has no edges between pins defined in a given $H_i$.

1. *Phase 1:*

   (a) i. Initially, the graph $H_1$ is identical to $G_1$.
   
   ii. The burnt graph $G_1^b$ for $G_1$ is constructed as follows: The pins in $G_1^b$ are the same as the pins in the original graph. There is an edge between two pins in $G_1^b$ iff there is an edge between the corresponding pins in $G_1$.
   
   (b) Repeat the following steps for $2 \leq i \leq n$.
   
   i. Let $A_i$ denote the set of all the terminals ( pins and explicit vertices). in $G_i$. Let the non-terminals called by $G_i$ be $G_{i_1} \cdots G_{i_k}$. Substitute the burnt graphs for each of the non-terminals called in $G_i$ to obtain $G_i'$. The cell $H_i$ is obtained as follows. The terminals in $H_i$ are in one-to-one correspondence with the terminals $A_i$ in $G_i$. Furthermore, there is an edge between two terminals iff either there was an edge between the corresponding terminals in the graph $G_i'$. $H_i$ also calls non-terminals $H_{i_1} \cdots H_{i_k}$ corresponding to the non-terminals $G_{i_1} \cdots G_{i_k}$ called in $G_i$. The one-one correspondence between the pins of non-terminals $H_{i_1} \cdots H_{i_k}$ and the terminal vertices of $H_i$ is the same as the one-one correspondence for $G_i$ except that for $G_{i_r}$, $1 \leq r \leq k$ we substitute $H_{i_r}$, $1 \leq r \leq k$.
   
   ii. Construct the burnt graph $G_i^b$ as follows: The pins in $G_i^b$ are the same as the pins in $G_i$. As in the case of $G_1^b$, there is an edge between two pins in $G_i^b$ iff there is an edge between the corresponding pins in $G_i'$.

2. *Phase 2:* Modify each $H_i$, $1 \leq i \leq n$ by removing any edges between pins in the definition of $H_i$.

3. Output $\Gamma_1 = (H_1, ..., H_n)$ as the new specification for $G$.

Figure 3: Algorithm for Producing Simple Specifications

builds on the well known vertex cover heuristic for the flat (non-hierarchical) case, where one computes a maximal matching and returns all the vertices involved in the matching as an approximate vertex cover. The algorithm in the non-hierarchical case has a performance guarantee of 2 [GJ79].

We note that the straightforward greedy approach for obtaining a maximal matching in a flat graph cannot be directly extended to the hierarchical case. Two reasons for this are as follows. First, the degree of a vertex in a hierarchical graph can be exponential in the size of the description, and so it is not possible to keep track of the neighbors of a node explicitly. Secondly, an edge between a pair of nodes can pass through several *pins*, and thus need not be explicitly present at any level. Therefore edges cannot be handled as simply as in the flat case. This complicates our heuristic since we can keep track of only a polynomial amount of information at each level.

Before we present the heuristic we give some notation which we use throughout this section. Given a graph $G$, $MM(G)$ denotes a maximal matching in the subgraph induced by the explicit vertices in $G$ (i.e. no pins and no nonterminals). $V(MM(G))$ denotes the vertices in the subgraph induced by $MM(G)$. $MxM(G)$ denotes a maximum matching of $G$ and $V(MxM(G))$ denotes the vertices in the subgraph induced by $MxM(G)$. We use $\psi(G_i)$ to denote the size of an approximate vertex cover for $E(G_i)$ (i.e. expanded version of $G_i$). We also use $EM(G_i)$ to denote the set of edges implicitly chosen by the heuristic from $E(G_i)$.

The following lemma recalls known properties of a maximum matching in a bipartite graph.

**Lemma 3.2** *Let $G = (S, T, E)$ be a bipartite graph and let $MxM(G)$ denote a maximum matching for $G$. Let $V_1^S$ and $V_1^T$ denote the set of vertices in $S$ and $T$ included in $V(MxM(G))$. Let $V_2^S$ and $V_2^T$ denote the set of vertices in $S$ and $T$ not included in $V(MxM(G))$. Then the following statements hold:*

1. *For all $\alpha \in V_2^S$ and $\beta \in V_2^T$, $(\alpha, \beta) \notin E$.*

2. *For all $v_x \in V_1^S$, $v_y \in V_1^T$, $v_z \in V_2^S$ and $v_w \in V_2^T$, if $(v_x, v_y) \in MxM(G)$ and $(v_y, v_z) \in E$ then $(v_x, v_w) \notin E$.*

**Proof:**

1. If $(\alpha, \beta) \in E$, then $\{(\alpha, \beta)\} \cup MxM(G)$ is also a feasible matching. This contradicts the assumption that $MxM(G)$ is a maximum matching for $G$.

2. Suppose $(v_x, v_y) \in MxM(G)$, $(v_y, v_z) \in E$, and $(v_x, v_w) \in E$. Then the matching $(MxM(G) - \{(v_x, v_y)\}) \cup \{(v_y, v_z), (v_x, v_w)\}$ contains more edges than $MxM(G)$, violating the assumption that $MxM(G)$ is a maximum matching. ■

Figure 4 gives the details of our approximation algorithm for minimum vertex cover.

**Heuristic HVC**
**Input:** *A simple hierarchical specification* $\Gamma = (G_1, ..., G_n)$ *of a graph* $G$.
**Output:** The size and a hierarchical description of an approximate vertex cover for $G$.

1. Repeat the following steps for $1 \leq i \leq n$.

   (a) Compute $MM(G_i)$.
   **Remark:** Recall that $MM(G_i)$ is a maximal matching on the subgraph of $G_i$ induced on the explicit vertices in $G_i$.

   (b) Compute $V_i^l$, where $V_i^l$ denotes the explicit vertices in $G_i$ not in $V(MM(G_i))$. Also let $G_i$ call non-terminals (if any) $G_{i_1}, ..., G_{i_k}$ in its definition. (Recall that $i_j < i$, $j = 1, 2, \cdots, k$.)
   **Remark:** Vertices in $V_i^l$ which are connected to pins in $G_{i_1}, ..., G_{i_k}$ are the endpoints of those edges that have their other endpoints in one of $G_{i_j}$ where $1 \leq i_j < i$.

   (c) **For** each vertex $v \in V_i^l$ **do**
   **If** $v$ is not adjacent to any nonterminals in $G_i$ **then** delete $v$ from $V_i^l$ **else**
   Let $v$ be adjacent to $p_{i_r} \in G_{i_r}^b$, such that $G_{i_r}$ is called in $G_i$.

      i. **If** there exists a marked edge incident on any of the $p_{i_r}$, $1 \leq r \leq k$, **then** match $v$ with $x_v$ such that $(x_v, p_{i_r})$ is a marked edge and delete $v$ from $V_i^l$ and $x_v$ from this copy of $G_{i_r}^b$.
      **else**

      ii. Choose a vertex $y_v$ such that $(y_v, p_{i_r})$ is an edge in $G_{i_r}^b$. Delete $v$ from $V_i^l$ and $y_v$ from this copy of $G_{i_r}^b$.

   (d) Let
   $$V_x^i = \{w \mid w \in V_i^l \text{ and } w \text{ is matched in step 1(c)}\}$$
   $$V_y^i = \{w \mid w \in V(G_{i_j}^b) \text{ and } w \text{ is matched in step 1(c)}\}$$

   (e) Construct a maximum matching on the set of vertices remaining in $V_i^l$ and the burnt graphs of nonterminals called in $G_i$.

   (f) For the bipartite graph $G_i^1$ induced by the vertices left over in $G_i$ including those in $G_{i_1}^b, ..., G_{i_k}^b$, and the pins in $G_i$, construct $MxM(G_i^1)$. $G_i^b$ for $G_i$ is the vertex induced subgraph of $MxM(G_i^1)$. The edges in $MxM(G_i^1)$ are **marked** in $G_i^b$.

   (g) $\psi(G_i) = |V(MM(G_i))| + |V_x^i| + |V_y^i| + \sum_{j=1}^{k} \psi(G_{i_j})$.
   **Remark:** Let $CM_i = \{(u,v) \mid u \in V_x^i, v \in V_y^i \text{ and } u \text{ and } v \text{ get matched up in Step 1(c)}\}$.
   $EM(G_i) = MM(G_i) \cup CM_i \cup \bigcup_{j=1}^{k} EM(G_{i_j})$. Note that $EM(G_i)$ is only needed in the proof; it is not explicitly computed. Further, $\psi(G_i) = 2 \times |EM(G_i)|$.

   (h) Construct $H_i$ as follows: The explicit vertices in $H_i$ are the vertices in the set $V(MM(G_i)) \cup V_x^i \cup V_y^i$. Their names are the same as those of the vertices in the sets $V(MM(G_i)) \cup V_x^i \cup V_y^i$. If $G_i$ calls a non-terminal $G_j$, $j < i$, then $H_i$ calls a copy of $H_j$.
   **Remark:** The $H_i$ created has the following property.
   Given a vertex $v \in E(G_i)$ as a path in the hierarchy tree, it is easy to check if $v$ occurs in $E(H_i)$ by simply following the same path. It is clear that if $v$ is in the approximate vertex cover then it will occur in a non-terminal on the path from the root to the non-terminal in which $v$ is defined.

2. Output $\psi(G_n)$ and the hierarchical specification $H = (H_1, ..., H_n)$.

Figure 4: Details of Vertex Cover Heuristic

## 3.3 Proof Of Correctness and Performance Guarantee

We now show that the above algorithm implicitly computes a maximal matching for $E(G_n)$.

**Lemma 3.3** $EM(G_n)$ *is a valid matching.*

**Proof:** We need to show that every vertex $u$ is in at most one edge in $EM(G_n)$.
**Case 1:** Vertex $u$ is matched with a vertex $v$ such that both $u$ and $v$ are explicitly defined in $G_i$, for some $i$, $1 \leq i < n$. This implies that in Step 1(b), the edge $(u,v)$ was chosen as an member of $MM(G_i)$. In Step 1(c) we do not consider any vertices which were in $V(MM(G_i))$. Hence $u$ is not an endpoint of any other edge in $EM(G_n)$.
**Case 2:** Vertex $u$ is matched with a vertex $v$ such that $u \in G_j$ and $v \in G_i$. Without loss of generality assume that $j < i$. In this case, $u$ was a part of the burnt graph $G_j^b$ and $G_i$ calls $G_j$. By Step 1(c), no edge incident on $u$ has been chosen in $MM(G_i)$. Once $(u,v)$ is chosen then in Step 1(c) we do not consider the vertices $u$ and $v$ anymore. ∎

**Lemma 3.4** *The matching $EM(G_n)$ is maximal.*

**Proof:** We need to prove that each edge in the expanded graph $E(\Gamma)$ has at least one of its endpoints in $EM(G_n)$. The proof consists of an exhaustive case analysis. Consider an edge $e \in E(T)$. There are two cases.
**Case 1:** Both endpoints of $e$ are explicit vertices in the definition of a cell $G_i$.

The proof for this case follows directly from Step 1 of the heuristic and the definition of $MM(G_i)$.
**Case 2:** Let $(v_i, v_j)$ denote the edge $e$ such that $v_i$ is in $G_i$ and $v_j$ is in $G_j$, where $j < i$. This edge $e$ passes through a sequence of pins $p_{i_r} \in G_{i_r}$, $1 \leq r \leq p$, where the path in the hierarchy tree from $G_i$ to $G_j$ consists of $G_{i_p}, \cdots G_{i_1}$ (see Figure 5). By the definition of hierarchical specification it is clear that for each pin in a nonterminal $G_k$ called in $G_t$, we have exactly one terminal in $G_t$ which is adjacent to the pin. We have two subcases to consider.
**Case 2.1:** $v_i \in V(MM(G_i))$ or $v_j \in V(MM(G_j))$. Here the proof follows from the definition of maximal matching.
**Case 2.2:** $v_i \notin V(MM(G_i))$ and $v_j \notin V(MM(G_j))$. We have two subcases again.
**Case 2.2.1:** $v_j \in V(MxM(G_j^b))$.

In this case we know that $v_j$ was matched with one of the pins. We have to consider two subcases depending on whether the vertex $v_j$ was a part of the burnt graph for all the non-terminal nodes on the path from $G_j$ to $G_i$ in the hierarchy tree, or it was a part of burnt graphs for some non-terminal and subsequently got dropped.
**Case 2.2.1.1:** $\forall m$ such that $1 \leq m \leq p$, $v_j \in V(G_{i_m}^b)$. (Informally, this means that the vertex $v_j$ was a part of the burnt graph for every non-terminal which is on the path from $G_i$ to $G_j$.)
In this case when we process the cell $G_i$ either $v_i$ or $v_j$ get matched up in Step 1(c). Hence the edge $(v_i, v_j)$ is covered.
**Case 2.2.1.2:** $\exists m$ $(1 \leq m < p)$ such that $v_j \in V(G_{i_{m-1}}^b)$ and $v_j \notin V(G_{i_m}^b)$. (Informally, $v_j$ was not part of the burnt graph for cell $G_{i_m}$, and $G_{i_m}$ is on the path from $G_i$ to $G_j$ in the hierarchy tree.)

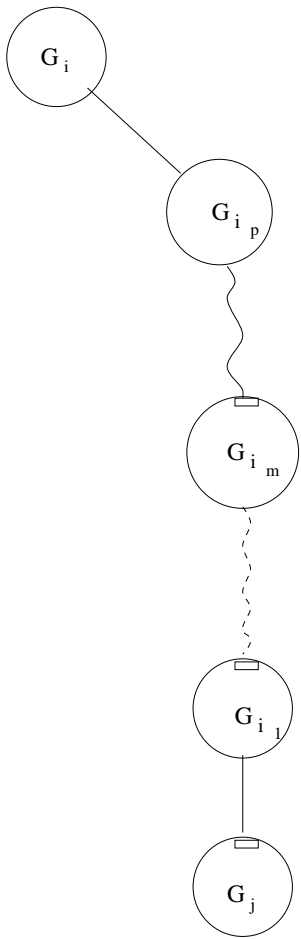

Figure 5: Figure showing the position of $G_i$ and $G_j$ in the hierarchy tree.

In this case, if $v_j$ gets matched with some other vertex, we are done. So, assume that $v_j$ is dropped (i.e. $v_j$ is not a part of the burnt graph). Now we need to show that $v_i$ gets a matching partner when it is picked up for processing. This case is more complicated and the proof uses the following lemmas (which, in turn, are proven using Lemma 3.2).

**Lemma 3.5** *Let $v_j$ be adjacent to pins $p_{i_1}^{i_m}, p_{i_2}^{i_m}, \cdots p_{i_k}^{i_m}$ in $G_{i_m}$ and let $v_j \notin V(G_{i_m}^b)$ (i.e. $v_j$ was not picked up as a matching partner for any of the pins). Then the following statements hold:*

1. *Each pin $p_{i_l}^{i_m}$ is matched with a distinct vertex $v_{i_l}^{i_m}$, $1 \leq l \leq k$.*

2. *$\forall l, 1 \leq l \leq k$, $v_{i_l}^{i_m}$ is not adjacent to any pin $p_r^{i_m}$ such that $p_r^{i_m}$ does not have a matching partner in $G_{i_m}^b$.*

**Proof:**
(1) Follows from the fact that we computed a maximum matching in Step 1(e) and (1) of Lemma 3.2.
(2) Follows from (2) of Lemma 3.2. ∎

Call a vertex $v_{i_l}^{i_m}$ a **private partner** of a pin $p_{i_l}^{i_m}$ in $G_{i_m}$, if $v_{i_l}^{i_m}$ is matched up with $p_{i_l}^{i_m}$ and is not adjacent to any pin $p_r^{i_m}$ in $G_{i_m}$ which does not have a matching partner. The following lemma says that if $v_j$ gets dropped off at stage $G_{i_m}$, each of the subsequent pins which are on the path from $v_j$ to $v_i$ has a private matching partner.

**Lemma 3.6** *Let $v_j \in V(G_{i_{m-1}}^b)$ and $v_j \notin V(G_{i_m}^b)$. Let $p_x^{i_m}$ be a pin in $G_{i_m}$, which is adjacent to $v_j$ and terminates at $v_i$. Then each of the pins $p_1^{i_m}(= p_x^{i_m}), p_2^{i_{m+1}}, \cdots, p_{p-m+1}^{i_p}$ on the path from $p_1^{i_m}$ to $v_i$ has a private partner in $G_{i_q}^b$, $m \leq q \leq p$.*

**Proof:** By induction on the length of the path from $G_i$ to $G_j$ in the hierarchy tree $HT(\Gamma_i)$.
**Basis:** The path is of length 1. By (1) and (2) of Lemma 3.5 it follows that $p_1^{i_m}$ has a private partner.
**Induction:** Assume that the Lemma holds for all paths of length $\lambda$. Now consider a path of length $\lambda + 1$. Again by Lemma 3.5, $p_1^{i_m}$ is matched up with say $v_k$. By (1) and (2) of Lemma 3.5, we know that $v_k$ is the private partner of $p_1^{i_m}$. We therefore have only two cases to consider.
**Case 1:** $v_k$ gets matched up with $p_2^{i_{m+1}}$.
In this case we can use our induction hypothesis and we are done.
**Case 2:** $v_k$ gets dropped.
By (1) and (2) of Lemma 3.5, we know that the pin $p_2^{i_{m+1}}$ will get some other private partner. Now, by Induction hypothesis we are done.□

We now continue the proof of Case 2.2.1.2. By Lemma 3.6 it follows that when $G_i$ is processed, pin $p_{p-m+1}^{i_p} \in G_{i_p}$ has a private partner. Therefore, when we process $v_i$, $v_i$ is sure to get matched up, because the private partner of $p_{p-m+1}^{i_p}$ which is adjacent to $v_i$ cannot be used as matching partner by any other vertex in $G_i$. So that the edge $(v_i, v_j)$ is covered by the vertex $v_i$.
**Case 2.2.2:** $v_j \notin V(MxM(G_j^b))$. The argument is similar to that of Case 2.2.1.2 because $v_j$ gets dropped at the very first stage. ∎

**Theorem 3.7** *Given a hierarchical graph $G$, the above approximation algorithm computes an approximate vertex cover within factor of 2 of the optimal value.*

**Proof:** Follows from Lemmas 3.3 and 3.4. ∎

## 3.4 Query Problem

We can easily modify our algorithm to answer the query problem. For this, we can use the hierarchical representation of the solution obtained.

**Theorem 3.8** *Given any vertex $v$ as a path in the hierarchy tree, we can determine in $O(N + M)$ if $v$ is in the approximate vertex cover so computed.*

**Proof:** Observe that the hierarchy tree for $H$ is identical to the hierarchy tree for $\Gamma$ except that the nodes in $HT(H)$ are labeled by $H_i$, whenever the corresponding node in $HT(\Gamma)$ is labeled $G_i$. This means, that the sequence of nonterminals used to identify the query vertex $v$ can be used to to check if $v$ is in the approximate vertex cover computed. For this, note that if $v$ is in the approximate vertex cover, then it lies on the path from the root of $HT(H)$ to a nonterminal $H_i$ such that $v$ is in the corresponding $G_i$ in the original graph $G$. ∎

The hierarchical specification can be used to output the approximate solution computed. For this, we do a simple preorder traversal of the nodes in the hierarchy tree $HT(H)$ and output the explicit nodes in each cell. Its easy to see that we can output the solution in $O(N)$ space (since the depth of $HT(H)$ no more than $n$ and each node on a path from root to a leaf is labeled with a distinct cell) and time linear in the size of $E(\Gamma)$.

## 3.5 Time Complexity

**Theorem 3.9** *HVC runs in time $O(N^{3.5})$.*

**Proof:** We compute a maximum matching at each level. It is well known that a maximum matching for a graph $G(V, E)$ can be found in time $O(|V|^{2.5})$ [MV80]. Thus computing a maximum matching while processing $G_i$ takes $O((n_i + \sum_{l=1}^{k} p_{i_l})^{2.5})$ time where $p_{i_1}, ..., p_{i_k}$ are respectively the number of pins in cells $G_{i_1}, ..., G_{i_k}$ which are called in the definition of $G_i$. We also compute a maximal matching while processing each $G_i$ and the time for this is $O(n_i + e_i)$, where $e_i$ is the number of edges in the level $i$. Therefore, the total time complexity is bounded by $\sum_{i=1}^{n}(O((n_i + \sum_{l=1}^{k} p_{i_l})^{2.5}) + O(n_i + e_i))$ which is bounded by $O(N^{3.5})$. ∎

**Corollary 3.10** *Given a hierarchical specification of a graph $G$, we can compute in time polynomial in the size of the specification, the size of an approximate minimum maximal matching which is within a factor of 2 of the optimal.*

**Proof:** Follows from the fact that any maximal matching is within a factor of 2 of the optimal minimum maximal matching. ∎

# 4 Approximating Weighted Max Cut

Given an undirected graph $G(V, E)$, the goal of the simple max cut problem is to partition the set $V$ into two sets $V_1$ and $V_2$ such that the number of edges in $E$ having one end point in $V_1$ and the other in $V_2$ is maximized [GJ79].

In [HR+93], it is shown that MAX CUT$_{HG}$ is PSPACE-hard. In this section, we show that given a hierarchical specification of a graph $G$, we can compute an approximate max cut which is within 2 times the optimum and a hierarchical specification of the vertices in one of the sets in the partition. Our algorithm computes the number of edges in the approximate cut in time polynomial in the size of the hierarchical description. An algorithm for weighted max cut can be devised along the same lines and is omitted. Since a graph obtained by expanding a hierarchical specification can in general be a multigraph, our approximation algorithms treat copies of an edge as distinct edges.

We begin with a brief overview of the algorithm. First, we recall the idea behind the known heuristic for computing a near optimal weighted max cut in the flat (non-hierarchical) case. That heuristic (referred to as FMAX-CUT in the following discussion) processes the nodes in arbitrary order, and assigns each node $v$ either to $V_1$ or to $V_2$ depending upon which of these sets has edges of least total cost to $v$. As in the case of the vertex cover algorithm, our approximation algorithm for MAX CUT$_{HG}$ processes the input specification in a bottom up fashion. At each level, we construct a burnt graph $G_i^b$ starting from the original description of the cell $G_i$. We use the heuristic FMAX-CUT to partition the explicit vertices at each stage. The burnt graph $G_i^b$ for $G_i$ then consists of two super nodes denoting an implicit partition of all the vertices defined in levels below. The edges go from a super node to the pins in $G_i$. Each edge has a weight associated with it. The edge weight is the number of edges the explicit vertex represented by the pin has to the vertices in that partition. In the following description, $A_i$ denotes the set consisting of all the explicit vertices in $G_i$ which are not adjacent to any nonterminals in the definition of $G_i$. Further, let $G(A_i)$ denote the subgraph induced on the nodes in $A_i$. The sets $V_1^i$ and $V_2^i$ denote the partition of the vertices of $E(G_i)$. Let $E^i$ denote the number of edges in the near optimal cut of $E(G_i)$. Also, for any vertex $v$, let $Count_v(V_j^i)$ denote the number of edges having one endpoint as $v$ and the other endpoint in the set $V_j^i$. Throughout this section, the reader should bear in mind that as a consequence of the definition of hierarchical specification, a terminal (an explicit vertex or a pin) defined in $G_i$ can be adjacent to at most one pin in each nonterminal called in $G_i$. The details of the approximation algorithm HMAX-CUT appear in Figure 6.

**Example:** Figure 7 illustrates the execution of the algorithm for the hierarchical specification given in Figure 1. The figure consists of 3 columns. The first column corresponds to $G_i$. The second column denotes the burnt graph $G_i^b$ of $G_i$. As mentioned before, the weights on the edges denote the number of vertices in $V_j^i$ that are adjacent to the pin. The third column shows the hierarchical representation $H$ being obtained level by level. ∎

## 4.1 Proof of Correctness

We now prove that the algorithm indeed produces a valid implicit partition of vertices.

**Theorem 4.1** *Given a hierarchical specification $\Gamma$, the heuristic HMAX-CUT computes a partition of the given vertex set.*

**Proof:** Induction on the number of non-terminals in the definition of $\Gamma$.
**Basis:** When $\Gamma = (G_1)$. In this case the theorem follows by the correctness of FMAX-CUT.

**Heuristic HMAX-CUT**
**Input:** *A simple hierarchical specification $\Gamma = (G_1, ..., G_n)$ of a graph $G$.*
**Output:** *A hierarchical specification $H = (H_1 \cdots H_n)$ of the vertices in the set $V_1^n$ and $E^n$ the number of edges in the approximate max cut computed.*

1. **For** $1 \leq i \leq n$ **do**

   (a) Use Algorithm FMAX-CUT to partition $A_i$ into sets $X_1^i$ and $X_2^i$. (Note that we do not consider any edges which are from these explicit vertices to the pins.)

   (b) $E_1^i$ = number of edges $(u, v)$ such that $u \in X_1^i$ and $v \in X_2^i$.

   (c) Let $G_i$ call nonterminals $G_{i_1}, \cdots, G_{i_m}$ in its definition. Let $B_i$ denote the set of all the explicit vertices remaining after Step 1(a). (Note that each of these explicit vertices is adjacent to at least one nonterminal in the definition of $G_i$.) We consider the vertices in $B_i$ one at a time. Let $Y_1^i = Y_2^i = \phi$.
   **For** each vertex $v \in B_i$ **do**

      i. Compute sets $V_{X_1^i}^v$ and $V_{X_2^i}^v$ defined by
      $V_{X_1^i}^v = \{w | w \in X_1^i \text{ and } w \text{ is adjacent to } v\}$ and $V_{X_2^i}^v = \{w | w \in X_2^i \text{ and } w \text{ is adjacent to } v\}$

      ii. **If** $G_i$ calls no nonterminals **then** $Count_v(V_1^i) = |V_{X_1^i}^v|$ and $Count_v(V_2^i) = |V_{X_2^i}^v|$ **else**
      Let $v$ be adjacent
      to pins $p_{v,i_l} \in G_{i_l}$, $1 \leq l \leq m$.
      Let $wt(V_1^{i_l}, p_{v,i_l})$ denote the weight of the edge between the super vertex $V_1^{i_l}$ and pin $p_{v,i_l}$. Then, let
      $$Count_v(V_1^i) = |V_{X_1^i}^v| + \sum_{1 \leq l \leq m} wt(V_1^{i_l}, p_{v,i_l})$$
      $$Count_v(V_2^i) = |V_{X_2^i}^v| + \sum_{1 \leq l \leq m} wt(V_2^{i_l}, p_{v,i_l})$$

      iii. **If** $(Count_v(V_1^i) \geq Count_v(V_2^i))$ **then** $Y_2^i = Y_2^i \cup \{v\}$ and $E_2^i = E_2 + Count_v(V_1^i)$
      **else** $Y_1^i = Y_1^i \cup \{v\}$ and $E_2^i = E_2 + Count_v(V_2^i)$

   (d) Construct the burnt graph $G_i^b$ as follows: The pins in $G_i^b$ are the same as the pins in $G_i$, and we have two super vertices $V_1^i$ and $V_2^i$ which implicitly represent the partition constructed so far. Let $G_i$ have $m_i$ pins in its definition. These pins will be connected to explicit vertices defined in $G_i$ and to pins in $G_{i_r}$, where $G_{i_r}$ is called in the definition of $G_i$. Let pin $p \in G_i$ be connected to pin $p_{i_r}$ in $G_{i_r}$. The weight of an edge $(p, V_j^i)$, $1 \leq r \leq m$, $1 \leq j \leq 2$, is calculated as follows:
   $$wt(p, V_j^i) = |Ex_j(G_i)| + \sum_{i_r} wt(p_{i_r}, V_j^{i_r})$$
   where $Ex_j(G_i) \subseteq X_j^i \cup Y_j^i$ denotes the set of explicit nodes in $G_i$ that are connected to $p$ and are added to $V_j^i$ in Steps 1(a) and 1(c).

   (e) $E^i = E_1^i + E_2^i + \sum_{i_r} E_{i_r}$

   (f) $H_i$ has no pins. The explicit vertices are in 1-1 correspondence with the vertices in the set $X_1^i \cup Y_1^i$. Furthermore, $H_i$ calls a non-terminal of type $H_{i_1} \cdots H_{i_m}$ corresponding to the nonterminals $G_{i_1} \cdots G_{i_m}$ called in $G_i$.
   **Remark:** Let $V_1^i = X_1^i \cup Y_1^i \cup \bigcup_{i_j} V_1^{i_j}$ and $V_2^i = X_2^i \cup Y_2^i \cup \bigcup_{i_j} V_2^{i_j}$, where $G_{i_j}$ $(i_j < i)$, $1 \leq j \leq m$ appears in the definition of $G_i$.

2. Output $E^n$ and the hierarchical specification $H = (H_1 \cdots H_n)$.

Figure 6: Details of MAX-CUT Heuristic

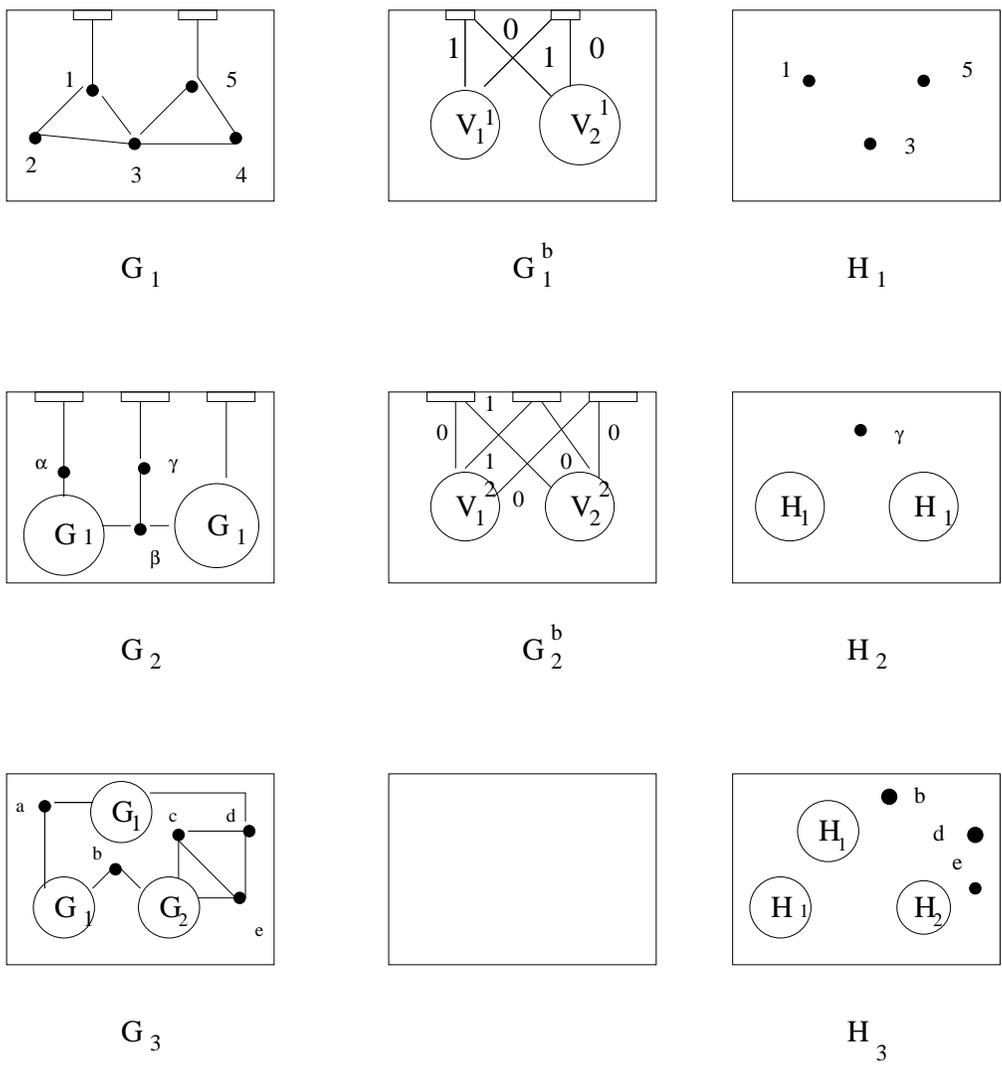

Figure 7: Figure showing the execution of heuristic HMAX-CUT on the specification given in Figure 1.

**Induction Step:** Assume that the theorem holds for all specifications with at most $(n-1)$ non-terminals. Consider the case when $\Gamma = (G_1, G_2, \cdots G_n)$. Let $G_n$ call non-terminals $G_{n_1}, G_{n_2}, \cdots, G_{n_k}$. By the induction hypothesis, we know that the the vertices in the hierarchy tree rooted at $G_{n_k}$ are partitioned into two sets. The explicit vertices of $G_n$ are clearly partitioned into two sets $X_1^n \cup Y_1^n$ and $X_2^n \cup Y_2^n$. Moreover, $V_r^n = X_r^n \cup Y_r^n \cup \bigcup_{n_j} V_r^{n_j}$, $1 \leq r \leq 2$. Therefore, it follows that the algorithm partitions the vertices into two sets. ∎

## 4.2 Performance Guarantee

We first prove that the weights on the edges in the burnt graph from super nodes to the pins actually represent the number of nodes in the definition that the pin is adjacent to.

**Lemma 4.2** *Let $\Gamma$ be a hierarchical specification of a graph $G$ constructed by HMAX-CUT. Consider the burnt graph $G_i^b$ corresponding to the non-terminal $G_i$ in the hierarchical specification. Then the weight of an edge from a pin $p \in G_i$ to a super vertex $V_j^i$, $1 \leq j \leq 2$, is equal to the total number of edges from $p$ to the vertices in the set represented by $V_j^i$.*

**Proof:** We prove the theorem for $V_1^i$. The proof for $V_2^i$ is similar. The proof is by induction on the number of nonterminals in the definition of $\Gamma$.

**Basis:** When $\Gamma = (G_1)$. In this case the lemma follows by fact that the weights were calculated by counting the number of explicit vertices in $G_1$ that are adjacent to the pin.

**Induction Step:** Assume that the lemma holds for all specifications which have no more than $(n-1)$ non-terminals. Consider the case when $\Gamma = (G_1, G_2, \cdots G_n)$. Let $G_n$ call $G_{i_1}, G_{i_2}, \cdots, G_{i_m}$. By the induction hypothesis, we know that the lemma holds for the burnt graphs corresponding to the non-terminals $G_{i_1}, G_{i_2}, \cdots, G_{i_m}$. Consider the non-terminal $G_n$. In Steps 1(a) and 1(c) the explicit vertices are partitioned into two sets $X_1^n \cup Y_1^n$ and $X_2^n \cup Y_2^n$. Consider a pin $p$ in $G_n$. Clearly, the total number of edges from $p$ to the vertices in the set $V_j^i$ is equal to $|Ex_1(G_n)| + \sum_{i_k} Edges_p(G_{i_k})$, where, $1 \leq k \leq m$ and $Ex_1(G_n) \subseteq X_1^n \cup Y_1^n$ represents the explicit vertices in $G_n$ that are adjacent to the pin $p$, and $Edges(G_{i_k})$ represents the number of edges which have one endpoint in $G_{i_k}$ and are incident on the pin $p$.

Note that the edges incident on the pin $p$ with one end point in $G_{i_k}$, $(1 \leq k \leq m)$ pass through the pins in the definition of $G_{i_k}$. By the induction hypothesis, the weight represents the number of edges from the pin to the explicit vertices defined in the graph $E(G_{i_k})$. The lemma now follows. ∎

We are now ready to prove that the heuristic computes a near-optimal maximum cut.

**Lemma 4.3** *Let $\Gamma$ be a hierarchical specification of a graph $G$. Let $\Xi_j$ denote the number of edges which are explicitly defined in $E(G_j)$. Then, $\Xi_n \leq 2E^n$*

**Proof:** The proof is by induction on the number of non-terminals in the hierarchical specification.

**Basis:** When there is only one non-terminal, the result follows by the correctness of the procedure FMAX-CUT.

**Induction Step:** Assume that the theorem holds for all hierarchical specifications which have no more than $(n-1)$ non-terminals in their definition. Consider the hierarchical specification $\Gamma = (G_1, G_2, \cdots G_n)$.

Consider the definition of the non-terminal $G_n$. Let $G_n$ call $G_{i_1}, G_{i_2}, \cdots, G_{i_k}$. The edges in $E(G_n)$ can be divided into three different categories.

1. Type 1 edges which have both the end points explicitly defined in one of the hierarchy trees rooted at $G_{i_r}$, $1 \leq r \leq k$.

2. Type 2 edges which have both the endpoints explicitly defined in the definition of $G_n$.

3. Type 3 edges which have one endpoint defined explicitly in $G_n$ and the other endpoint defined in a non-terminal occurring in one of hierarchy tree $HT(G_{i_r})$ rooted at $G_{i_r}$, $1 \leq r \leq k$.

Also let $Exp_j$ denote the number of edges which occur explicitly in the definition of $G_j$. Then clearly the total number of edges $\Xi_n$ equals,

$$\Xi_n = \sum_{i_r} \Xi_{i_r} + Exp_n + Cross_n$$

where $Cross_n$ denotes the set of Type 3 edges. By induction hypothesis, we know that the vertices in the hierarchy tree rooted at $G_{i_k}$ are partitioned into two sets such that the number of edges crossing the cut is at least 1/2 of the total number of edges. Therefore, $\forall i_r, \Xi_{i_r} \leq 2E_{i_r}$. By Step 1(c), explicit vertices in $G_n$ which are not adjacent to any pins are partitioned in such a way that the at least half of the of edges in the subgraph induced by these vertices are cut. Each remaining explicit vertex in $G_n$ is added to the set $V_1^n$ or $V_2^n$ depending on which set has fewer vertices adjacent to it. By Lemma 4.2, the weights on the edges from the pins to the super nodes by give the number of nodes that the pin is adjacent to in the hierarchy tree rooted at that non-terminal. Therefore, $Exp_n + Cross_n \leq 2(E_1^n + E_2^n)$, and hence

$$\Xi_n = \sum_{i_r} \Xi_{i_r} + Exp_n + Cross_n \leq 2E^n$$

. ∎

**Theorem 4.4** *Let $\Gamma$ be a hierarchical specification of a graph $G$. Let $OPT(G)$ denote a maximum cut in $E(G)$. Then $|OPT(G)| \leq 2E^n$.*

**Proof:** The theorem follows from the above lemma and the fact that $|OPT(G)| \leq \Xi_n$. ∎

### 4.3 Query Problem

Using the above hierarchical specification of the set of vertices in $V_1^n$, we can answer the question of which set a given vertex belongs. As mentioned earlier, we assume that a vertex $v$ is specified as a sequence of nonterminals which occur on the path from the root to the nonterminal in which $v$ occurs.

**Theorem 4.5** *Let $\Gamma$ be a hierarchical specification of a graph $G$ with $N$ vertices. Given any vertex $v$ in the graph $G$, we can determine in $O(N)$ time, the set to which $v$ belongs.*

**Proof:** Observe that, the hierarchy tree $HT(H)$ of $H$ is identical to $HT(\Gamma)$ except that if a node in $HT(\Gamma)$ is labeled by $G_i$ then the corresponding node in $HT(H)$ is labeled by $H_i$. This means that the sequence of nonterminals used to specify $v$ in $E(\Gamma)$ can be directly used to locate the nonterminal $H_i$ in which $v$ may occur. This implies that, given a vertex $v$ one can easily check in $O(N)$ time if the vertex occurs in $H$ by following the path in the hierarchy tree to the non-terminal in which $v$ occurs. If $v$ appears in $H$ then it belongs to the set $V_1^n$, else it is in the set $V_2^n$. ∎

As in the case of vertex cover problem, the hierarchical specification $H$ obtained can be used to output the $V_1^n$. For this, we do a simple preorder traversal of the nodes in the hierarchy tree $HT(H)$ and output the explicit nodes in each cell. This takes $O(N)$ space and time linear in the size of $E(\Gamma)$.

### 4.4 Time Complexity

**Theorem 4.6** *The algorithm HMAX-CUT runs in time $O(N + M)$ and constructs a hierarchical specification of size $O(N)$ of the set $V_1^n$.*

**Proof:** Consider the time taken to process $G_i$. Step 1 (a) takes $O(n_i + m_i)$ time. Steps 1 (c) and 1 (d) take $O(d_j) + O(1)$ time to process each terminal of degree $d_j$ in $G_i$. Therefore, the total running time of Steps 1 (c) and (d) is $O(n_i + m_i)$. Hence the total running time of the algorithm is $\sum_{1 \leq i \leq n} O(n_i + m_i) = O(N + M)$. Size of each $H_i$ is no more than $n_i$, the number of vertices in $G_i$. Hence the size of $H$ is $\sum_i O(n_i) = O(N)$. ∎

## 5 Approximating Bounded Degree Maximum Independent Set

Our heuristic for obtaining a near-optimal solution to the maximum independent set problem on bounded degree hierarchically specified graphs is based on a well known heuristic in the flat case. The heuristic in the flat case (referred to FIND-SET in the subsequent discussion) is the following. We pick and add an arbitrary node $v$ to the approximate independent set and delete $v$ and all the nodes which are adjacent to $v$. This step is repeated until no nodes are left. It is easy to see that for a graph in which each node has degree at most $B$, the independent set produced by the heuristic is within a factor $B$ of the optimal value. We now show how to extend this heuristic to the hierarchical case. Throughout this section, we use $V_j$ to denote the set of vertices from $E(G_j)$ that are in the approximate independent set produced by the algorithm. The details of the heuristic HIND-SET are given in Figure 8.

### 5.1 Performance Guarantee and Proof Of Correctness

We now show that the approximate independent set computed is within a factor of $B$ of the optimal independent set.

**Lemma 5.1** *The set $V_n$ produced by HIND-SET is a maximal independent set.*

**Proof:** The proof follows by an easy induction on the number of non-terminals in the hierarchical specification $\Gamma$. ∎

**Lemma 5.2** *Let $OPT(G)$ denote the size of an optimal independent set in $G(= E(\Gamma))$. Then $|V_n| \geq \frac{OPT(G)}{B}$.*

**Proof:** Follows from the fact that every time we choose a vertex, we delete (mark) no more than $B$ terminals (explicit vertices and pins). ∎

**Heuristic HIND-SET**
**Input:** A simple hierarchical specification $\Gamma = (G_1, ..., G_n)$ of a graph $G$. Each node of $G$ has a degree of at most $B$, where $B$ is a constant.
**Output:** A *hierarchical specification* $H = (H_1, ..., H_n)$ of the approximate independent set and $|V_n|$, the size of the approximate independent set.

1. Repeat the following steps for $1 \leq i \leq n$.

   (a) Let $A_i$ denote the set of all the explicit vertices in $G_i$. Starting from the set $A_i$, we create a new set $B_i$ as follows. For each vertex $v \in A_i$, we place it in the set $B_i$ iff $v$ is *not* adjacent to any of the pins marked **removed** in the burnt graphs of $G_j$, where $G_j$, $j < i$, appears in the definition of $G_i$. Let $G(B_i)$ denote the subgraph induced on the nodes in $B_i$.

   **Remark:** A vertex $v$ is placed in the set $B_i$ iff none of its neighbors in $G_j$, $j < i$, have been placed in $V_j$.

   (b) Use Algorithm FIND-SET on $G(B_i)$ to obtain the independent set $X_i$.

   **Remark:** We do not consider any edges which are from these explicit vertices to the pins.

   (c) Let $|V_i| = |X_i| + \sum_j |V_j|$ where $G_j$, $j < i$, appears in the definition of $G_i$.

   **Remark:** $V_i = X_i \cup \bigcup_j V_j$ where $G_j$, $j < i$ appears in the definition of $G_i$. (Observe that the set is created implicitly.)

   (d) Now construct the burnt graph $G_i^b$ for $G_i$ as follows: The pins in $G_i^b$ are the same as the pins in $G_i$. A pin in $G_i$ is marked **removed** iff the pin is *either* adjacent to one of vertices in the set $X_i$ *or* it is adjacent to one of the pins in $G_j$ ($j < i$), which is marked **removed**.

   (e) Construct $H_i$ as follows: The explicit vertices in $H_i$ are the vertices in the set $X_i$. If $G_i$ calls a non-terminal $G_j$, $j < i$, then $H_i$ calls $H_j$.

2. Output $|V_n|$ as the size of approximate independent set and $H = (H_1, ..., H_n)$ as the hierarchical specification of the approximate independent set.

Figure 8: Details of Heuristic for Maximum Independent Set

## 5.2 Query Problem

As in the case of max cut problem, the hierarchy tree of $H$ is identical to the hierarchy tree $HT(\Gamma)$ of $\Gamma$, except that the corresponding nodes are labeled by $H_i$ instead of $G_i$.

**Theorem 5.3** *Let $\Gamma$ be a hierarchical specification of a graph $G$. Given any vertex $v$ in the graph $G$, we can determine in $O(N)$ time, if $v$ belongs to the approximate independent set obtained.*

**Proof:** Given the label of any node as a path in the hierarchy tree, it is easy to check if the vertex belongs to the independent set specified by $H$. This can be done by traversing the hierarchy tree $HT(H)$ and checking if the vertex appears in the given $H_i$. ∎

As in the case of previous algorithms, we can output the solution in $O(N)$ space and time linear in the size of $E(\Gamma)$. This can be done by a preorder traversal of the hierarchy tree $HT(H)$.

## 5.3 Time Complexity

**Lemma 5.4** *The algorithm HIND-SET runs in time $O(N + M)$ and constructs an $O(N)$ size hierarchical specification for the approximate independent set.*

**Proof:** The proof follows by observing that HIND-SET processes each of the $G_i$ in $O(n_i + m_i)$ time. ∎

Summarizing the above results, we have:

**Theorem 5.5** *Let $\Gamma$ be a hierarchical specification of a graph $G$ with maximum node degree $B$. Then we can compute in time $O(N + M)$ (the size of the specification), an approximate independent set which is within a factor $B$ of the size of a maximum independent set.* ∎

## 6 Approximating Weighted MAX 3SAT

We now consider the problem of finding a truth assignment to the variables of a hierarchically specified instance of 3SAT so as to maximize the number of clauses that can be simultaneously set to true. We first outline a heuristic (see Figure 9) with performance guarantee 2, which works for non-hierarchical specifications of MAX 3SAT instances. The heuristic is a variant of a heuristic for MAX 3SAT in [Jo74].

We first observe that the approximation algorithm given in Figure 9 has a performance guarantee of 2.

**Lemma 6.1** *Let $|C|$ denote the number of clauses in $F$. Let $Heu(F)$ denote the number of clauses set true by FMAX 3SAT. Then $Heu(F) \geq |C|/2$.*

**Proof:** Let $C_{x_i}$ denote the number of clauses in the star centered around $x_i$. We know that the value assigned to $x_i$ in Step 2(b) satisfies at least $C_{x_i}/2$ clauses. Given that $\sum_{x_i} C_{x_i} = |C|$, the lemma follows. ∎

Next we show how, given a hierarchical specification of a 3SAT formula $f$ we can construct a hierarchical specification of the bipartite graph corresponding to $f$. The transformation is given in Figure ??.

It is easy to see that the transformation given in Figure ?? constructs a hierarchical specification of the bipartite graph associated with the 3SAT formula $f$. Thus we have:

**Lemma 6.2** *Given an instance $F = (F_1(X^1), \ldots, F_{n-1}(X^{n-1}), F_n)$ of $3SAT_{HG}$. Procedure TFORM constructs a hierarchical specification $BG(F) = (G_1, \ldots, G_n)$ such that*

1. *size of $BG(F)$ is $O(size(F))$.*

2. *$BG(F)$ can be constructed in $O(size(F))$ time.*

3. *$E(BG(F))$ is the bipartite graph associated with the formula $E(F)$.* ∎

The basic idea of the approximation algorithm for the hierarchical case is to mimic the flat case algorithm FMAX 3SAT. The approximation algorithm is fairly simple, and its details appear in Figure 9.

In the rest of the section, we let $A_i$ be the set consisting of all variables in $F_i$ which are not adjacent to any nonterminals in the definition of $F_i$. Further, let $F(A_i)$ denote the subgraph induced on the nodes in $A_i$. The details of the heuristic HMAX-3SAT appear in Figure ??.

## 6.1 Proof of Correctness and Performance Guarantee

The proof of the fact that the above algorithm guarantees a solution which is within 2 of the optimal value is easy and follows by verifying the following two lemmas which can easily be proven by an induction on the number of nonterminals in the definition of $\Gamma$.

**Lemma 6.3** *Each variable in the 3SAT formula $F$ specified by $\Gamma$ is assigned a unique truth value.* ∎

**Lemma 6.4** *Let $\Gamma = (F_1, F_2, \cdots, F_n)$ be a hierarchical specification of a 3SAT formula $F$. Consider the burnt graph corresponding to a non-terminal $F_i$ in the hierarchical specification. Then the weight of an* *edge from a pin $p_i$ to the super vertex $P_i$ ($N_i$) represents the total number clauses in which the variable represented by $p_i$ occurs un-negated (negated) in the expanded formula denoted by $E(F_i)$.* ∎

By an easy induction on the number of nonterminals in the definition of $\Gamma$ and using the above lemmas we can prove that

**Theorem 6.5** *Heuristic HMAX 3SAT has a performance guarantee of 2.* ∎

## 6.2 Query Problem

We show that the algorithm given above can in fact be used to give a hierarchical description of the truth assignments to the variables of the 3SAT formula $F$.

**Theorem 6.6** *Let $\Gamma$ be a hierarchical specification of a 3SAT formula $F$. Given a variable $v$ in the 3SAT formula we can tell in $O(N)$ time, the truth assignment to the variable $v$.*

**Proof:** To do this we simply follow the path from the root to the nonterminal in which the variable occurs, and then check the truth value assigned to $v$. Since each clause has at most 3 variables, we can also tell the truth value of any clause in $F$ in $O(N)$ time. ∎

## 6.3 Time Complexity

**Theorem 6.7** *Given a hierarchical specification of a 3SAT formula $f$, the algorithm HMAX-3SAT runs in time $\sum_{1 \leq i \leq n} O(n_i + m_i)$ and constructs a hierarchical specification of size $\sum_{1 \leq i \leq n} O(n_i)$ of the satisfying assignment to the variables in $f$, such that at least 1/2 total number of clauses in $E(\Gamma)$ are satisfied.*

**Proof:** Consider the time to process a cell $F_i$. If a vertex corresponding to a variable $v_j$ has degree $d_j$ in the definition of $F_i$, then it takes $O(d_j)$ time to find a truth assignment to $v_j$. Therefore, the total running time of Steps 1 (b) and (c) is $O(n_i + m_i)$. Hence the total running time of the algorithm is $\sum_{1 \leq i \leq n} O(n_i + m_i)$. Size of each $H_i$ is $n_i$, the number of vertices in $G_i$. Hence the size of $H$ is $\sum_i O(n_i)$. ∎

# 7 Non-Approximability Results

In this section we discuss our results on the non-approximability of several natural problems studied in the literature, when instances are specified hierarchically. We show that approximating the number of true gates in a hierarchically specified monotone acyclic circuit is PSPACE-hard. We then show that unless P = PSPACE the optimization versions of the high degree subgraph problem and the high vertex and edge connectivity problems cannot be approximated to within a factor $c < 2$.

Intuitively, problems proven to be P-hard by a *local* reduction (i.e. a reduction where each gate is replaced by a corresponding subgraph or gadget of *fixed size*), by a log-space reduction from MCVP, can be shown to PSPACE-hard by a polynomial time reduction from $\text{MCVP}_{HG}$. Such a reduction, transforms the given hierarchical specification of a monotone acyclic circuit **level by level** to obtain a hierarchical specification of the original problem instance. The proofs for the non-approximability of the optimization versions of the circuit value problem, high degree subgraph problem and the high-vertex and edge connectivity problems in the non-hierarchical case are examples of such *local* reductions from MCVP. This property of local reduction allows us to *lift* these reductions to the case when the inputs are specified hierarchically.

## 7.1 Approximating Number of True Gates in MVCP

The Monotone circuit value problem is known to be PSPACE-hard when the circuit is specified hierarchically [LW92, RH93]. We first observe that the problem is PSPACE-hard even for strongly 1-level-restricted hierarchical specifications.

**Lemma 7.1** *The problem MCVP is PSPACE hard even for strongly 1-level-restricted specifications in which a non-terminal $C_i$ calls exactly 2 copies of $C_{i-1}$.*

**Proof:** Follows from the fact that the instance of MCVP obtained by [LW92] in their reduction from QBF is of the required form. ∎

Before we give the PSPACE-hardness proof for $\text{MTG}_{HG}$, it is instructive to recall the proof by Serna [Se91], showing that MTG is P-complete. The proof consists of a log-space reduction from MCVP. Given an instance $C$ of MCVP with $n$ gates, the instance $C'$ of MTG consists of the same circuit $C$ along with $\lceil \frac{n}{\epsilon} \rceil$ additional AND gates forming a chain, with the first element of the chain being connected to the output gate of $C$ and the last element of the chain serving as the output for $C'$. As the circuit added to $C$ only propagates the value of output of $C$ it follows that

1. If $C$ outputs 0, then $OPT(MTG) < n$;

2. If $C$ outputs 1, then $OPT(MTG) \geq \lceil \frac{n}{\epsilon} \rceil$.

It is clear that the reduction can be done in log-space. As discussed in [Se91], the result holds even when instances are restricted to be planar.

We extend this result and show that $\text{MTG}_{HG}$ cannot be approximated to within any exponential function of the optimal. To show this, the basic idea is to construct a a chain of exponential number of AND gates using a simple specification, and join this chain in series to the output of an instance of $\text{MCVP}_{HG}$.

**Theorem 7.2** *Unless P=PSPACE, no polynomial time algorithm can approximate the maximum number of true gates in $MCVP_{HG}$ to within any $\eta^\epsilon$ ($\epsilon > 0$) factor of the optimal, even for simple strongly 1-level restricted hierarchical specifications, where $\eta$ denotes the size of the hierarchical specification.*

**Proof:** Let $C = \{C_1, C_2, ..., C_n\}$ be an instance of a simple hierarchical specification of $MCVP_{HG}$ in which each $C_i$ calls exactly two copies of $C_{i-1}$. Let $m$ denote the number of gates in $C$ and $N$ denote the size of $C$. We construct an instance $D = \{D_1, D_2, ..., D_n\}$ of a simple hierarchical specification of $MVCP_{HG}$ with $m + 2^{N^2}$ gates such that,

1. If $C$ outputs 0, then $OPT(D) < 2^N$;

2. If $C$ outputs 1, then $OPT(D) \geq 2^{cN^2}$, for some $0 < c \leq 1$.

We now discuss the construction of the instance $D$.

**Circuit $D_1$:** The Circuit $D_1$ consists of two disjoint circuits $D_{1,1}$ and $D_{1,2}$. $D_{1,2}$ is identical to $C_1$. The circuit $D_{1,1}$ has AND gates connected in series. The input of the first AND gate is connected to two pins. Similarly, the output of the last AND gate is connected to two pins. $D_{1,3}$ consists of a series of AND gates such that the total number of AND gates in $D_{1,1}$, $D_{1,2}$, $D_{1,3}$ equals $N \cdot n_1$. Figure 12 gives a schematic of the above construction.

**Circuit $D_i$, $2 \leq i \leq n-1$:** The Circuit $D_i$ consists of five circuits $D_{i,1}$, $D_{i,2}$, $D_{i,3}$, $D_{i,4}$ and $D_{i,5}$. $D_{i,4}$ and $D_{i,5}$ are identical to $D_{i-1}$. The circuits $D_{i,1}$ and $D_{i,2}$ each consists of a single AND gate. The AND gate corresponding to $D_{i,1}$ gets its input from two pins and its output is connected to the partial chain of AND gates in $D_{i,4}$. The AND gate corresponding to $D_{i,2}$ gets its input from the partial chain of AND gates in $D_{i,5}$ and its output is connected to a set of pins. $D_{i,3}$ consists of a series of AND gates and joins the partial chains of AND gates in the two copies of $D_{i-1}$. The total number of AND gates in $D_{i,1}$, $D_{i,2}$, $D_{i,3}$ equals $N \cdot n_i$. Figure 12 shows the schematic diagram of $D_i$.

**Construction of $D_n$:** As in $D_{n-1}$, $D_n$ consists of five circuits $D_{n,1}$, $D_{n,2}$, $D_{n,3}$, $D_{n,4}$ and $D_{n,5}$. $D_{n,4}$ and $D_{n,5}$ are identical to $D_{n-1}$. $D_{n,3}$ consists of a series of AND gates and joins the partial chains of AND gates in the two copies of $D_{n-1}$. The circuits $D_{n,1}$ and $D_{n,2}$ each consists of a single AND gate. The input port of the AND gate corresponding to $D_{n,1}$ is joined to the output port of $C$ and the output port feeds into the partial chain of the AND gates in $D_{n,4}$. The output of the AND gate corresponding to $D_{n,3}$ is designated as the output of $D$, and the input ports of $D_{n,3}$ are joined to the partial chain of AND gates in $D_{n,5}$. $D_{n,3}$ consists of a series of AND gates such that the total number of AND gates in $D_{n,1}$, $D_{n,2}$, $D_{n,3}$ equals $N \cdot n_n$. The construction is depicted in Figure ??.

Note that the size of $D$ denoted by $\eta$ is $O(N^2)$. Now, observe that the above construction specifies a circuit in which the output of the circuit corresponding to $C$ is connected to a exponentially long chain of AND gates. Given this observation it is not difficult to verify that the following lemma holds:

**Lemma 7.3** *If the output of $C$ is 1, at least $2^{cN^2}$ AND gates will output a 1; otherwise, less than $2^N$ of those gates will output a 1.*

Given Lemma 7.3 and the fact that the above construction of $D$ can be done in polynomial time the theorem follows. ∎

## 7.2 Approximating the Objective Function of a Linear Program

We now discuss our result concerning the nonapproximability optimizing the objective function of a hierarchically specified linear program. The PSPACE-hardness proof consists of it lifting the proof in [Se91] showing that approximating the objective function of a linear program is log-complete for $P$.

**Theorem 7.4** *Unless P=PSPACE, no polynomial time algorithm can approximate the objective function of an HLP to within any $\eta^\epsilon$ of the optimum, even for strongly 1-level restricted simple specifications. Here $\eta$ denotes the size of the specification.*

**Proof:** The reduction is from an instance of strongly 1-level-restricted simple hierarchical specification $D = \{D_1, D_2, ..., D_n\}$ of the problem $\text{MTG}_{HG}$. We construct an instance of $\text{LP}_{HG}$ $F = \{F_1, F_2, ..., F_n\}$, bottom up level by level as follows.

**Construction of $F_i$, $1 \le i \le n$:** Recall that the formula $F_i$ is of the form

$$F_i(X^i) = (\bigcup_{1 \le i_j \le i} F_{i_j}(X_j^i, Z_j^i)) \bigcup f_i(X^i, Z^i)$$

.

$$\Delta_i = \sum_{i_j} d_{i_j} \cdot \Delta_{i_j} + \sum_{z_j \in Z^i} c_j \cdot z_j$$

where, $\Delta_i$ is the objective function. We now describe each of the components in the above definition of $F_i$.

1. The set of dummy variables $X^i$ is in 1-1 correspondence with the pins of $F_i$. (Note that this implies that $X^n = \phi$.)

2. $Z^i = A^i \cup B^i$, where

   - $A^i = \cup_{i_r} A_{i_r}$ where the variables in $A_{i_r}$ are in 1-1 correspondence with the edges incident on the non-terminal $D_{i_r}$ called in $D_i$.

- The set $B^i$ consists of variables which are in 1-1 correspondence with the explicitly defined gates in $D_i$ and the 0 1 input ports of the circuit.

3. For the function $\Delta_i$ the coefficients are $c_i$ and $d_i$ are all 1.

4. $\forall i_r \ X_i^{i_r} = \phi$. (Note: This is true because the given circuit specification is 1-level-restricted.)

5. Corresponding to each $D_{i_r}$, $i_r < i$, called in $D_i$, we have a call to $F_{i_r}$ and the set of variables $Z_i^{i_r} \subseteq Z^i$ passed to $F_{i_r}$ are in 1-1 correspondence with the set of explicit variables in $F_i$ which correspond to the explicit gates defined in $D_i$.

We now describe the set of inequalities corresponding to $f_i(X^i, Z^i)$. We have one set of inequalities for each explicit gate in $F_i$. We also have an additional set of inequalities with each pin that is connected to the output port of an explicit gate in $D_i$. (The inequalities are very similar to those given in [Se91].)

1. If $x_k$ corresponds to an input port of the circuit, then we have the equation $x_k = 1$ if the corresponding input is 1 and the equation $x_k = 0$ if the corresponding input is 0.

2. For an AND gate, we have the inequalities $x_k \leq x_j$, $x_k \leq x_i$, $x_k \geq x_i + x_j - 1$, where $x_k$ is the variable denoting the AND gate and $x_i, x_j$ are the variables corresponding to the gates whose outputs serve as the inputs for the AND gate. If the gate is connected to a nonterminal, the variables $x_i$ and $x_j$ correspond to the variables that are associated with the edge joining the gate to the nonterminal.

3. For an OR gate, we have the inequalities $x_i \leq x_k$, $x_j \leq x_k$, $x_k \leq x_i + x_j$, where $x_k$ is the variable denoting the OR gate and $x_i, x_j$ are the variables corresponding to the gates whose outputs serve as the inputs to the OR gate.

4. Recall that with each pin we have an associated dummy variable. Consider a pin $p_j^i$ whose associated dummy variable is $x_j^i$. If $p_j^i$ is connected to the output port of a gate $x_k$ then we generate the equation $x_k = x_j^i$.

5. For each variable $x_k$ which denotes an edge going from an explicit gate to a nonterminal (i.e. $x_k$ is a variable in the set $A^i$) and is connected to an output port of an explicit gate, we generate the equation $x_k = x_j$ where $x_j$ denotes the variable corresponding to the gate which has an edge corresponding to $x_k$ joined to a nonterminal.

It is easy to see that the reduction gives rise to a *simple* strongly 1-level restricted specification of $F$, given that $D$ was simple and strongly 1-level restricted. Also, it is easy to see that the reduction can be done in polynomial time. Next observe that the reduction gives rise to a hierarchical specification $F$ which represents the set of inequalities which would be produced if the specification is expanded and Serna's construction [Se91] applied on the expanded circuit. The only difference that we have some intermediate variables on edges. Let $N$ be the size of $D$. The size of $F$, denoted by $\eta$, is $O(N^2)$.

Given the above observations, it is easy to verify that the value of $\Delta$ is less than $2^{2N}$ if the output of the circuit is 0 and the value of $\Delta$ is at least $2^{cN^2}$ for some $0 < c \leq 1$ if the output of the circuit is 1. The theorem follows. ∎

**Example:** Consider the hierarchical specification $D$ as given in Figure ??. The corresponding specification $F$ is given as follows:

$$\begin{aligned} F_1(x_1, x_2, x_3, x_4) &= \{(z_1 = x_1 \wedge x_2), (z_1 = x_3 = x_4)\} \\ F_2 &= \{(z_2 = 0), (z_3 = 1), (z_4 = z_2 \wedge z_3)\} \\ &\quad \cup F_1(a, b, c, d) \cup F_1(e, f, g, h) \cup \\ &\quad \{(z_4 = a = b), (z_5 = c \vee d)\} \cup \\ &\quad \{(z_5 = e = f), (z_6 = g \wedge h)\} \end{aligned}$$

Note that each equation involving an AND or an OR operator has to be replaced by the set of inequalities as discussed earlier.

The corresponding $\Delta$ function is also created similarly and is just a sum of all the explicit variables. Observe that the specification obtained is strongly 1-level-restricted and simple.

## 7.3 Approximating Connectivity and High Degree Subgraph Problems

Next, we consider the problems $\kappa$-HVCP, $\kappa$-HECP, and $k$-HDSP, when instances are specified hierarchically. We prove PSPACE-hardness results for these problems when instances specified hierarchically by *lifting* the known proofs showing the P-hardness of the corresponding problems in the non-hierarchical case. We illustrate this idea by presenting the PSPACE-hardness proof for $\kappa$-HVCP. PSPACE-hardness proofs for the other two problems are along the same lines.

The proof given in [KSS89] showing that $\kappa$-HVCP is P-complete is a log-space reduction from MCVP with additional restriction that outdegrees of all gates and the input nodes is at most 2, and there is at least one input node with whose value is 1. It can be easily shown by slightly modifying the reduction in [LW92] that

**Lemma 7.5** *The problem $MCVP_{HG}$ is PSPACE-hard even for hierarchical specifications satisfying all the following restrictions.*

1. *The specification is simple.*

2. *The specification is strongly 1-level-restricted.*

3. *Each $C_i$ calls exactly two copies of $C_{i-1}$.*

4. *The outdegree of all gates and the input nodes is at most 2.*

5. *There is at least one input node with whose value is 1.*

6. *The inputs and the outputs all occur in the last cell.* ∎

We recall the construction from [KSS89] to show the P-completeness of the $\kappa$-HVCP problem. Given an instance $C$ of the MCVP with the restriction that the outdegree of all gates and the input nodes is 2 and there is at least one input node with whose value is 1, an instance of $G$ 3-HVCP is created as follows:

1. Each input node of the circuit as well as the output node is replaced by a $K_{2,2}$ graph, as depicted in Figure ??(a).

2. Each OR gate of $C$ is replaced by a copy of the graph depicted in Figure ??(e). The upper nodes are called the *in-nodes* and the lower ones are referred to as the *out-nodes*.

3. Each AND gate of $C$ is replaced by a copy of the graph depicted in Figure ??(d).

4. An additional node $v_{new}$ is added and is connected to the out-nodes of the subgraph used to replace the output gate and all the in-nodes of the subgraphs replacing the input gates with value 1. The construction is illustrated through an example in Figure 16.

Using this construction it can be proven (see [KSS89]) that the output of $C$ is 1 iff the $G$ contains a 3-connected subgraph. As in the previous proof of PSPACE-hardness, we lift the reduction in the non-hierarchical case, to prove the PSPACE-hardness of 3-HVCP$_{HG}$.

**Theorem 7.6** *The problem $\kappa$-HVCP$_{HG}$ is PSPACE-hard for simple strongly 1-level-restricted hierarchical specifications.*

**Proof:** We prove the theorem for $\kappa = 3$. Given an instance $C = \{C_1, C_2, ..., C_k\}$ of simple hierarchical specification of $MCVP_{HG}$ in which each $C_i$ calls exactly two copies of $C_{i-1}$, we construct a simple hierarchical specification $\Gamma = \{G_1, G_2, ..., G_n\}$ of a graph $G$ such that $G$ has a 3-connected subgraph iff the circuit corresponding to $C$ outputs a 1. The reduction follows the same outline as in the proof of Theorem 7.2. It is done level by level and at each stage the gates of the circuit are replaced by a gadget depending on whether it is an AND or an OR gate.

**Graph $G_1$:** Except for a minor modification, the graph $G_1$ is the same as the one obtained using the construction (given above) proving the P-completeness of the problem in the flat (non-hierarchical) case. The modification is that if a gate in $C_1$ has its inputs connected to pins then the corresponding in-nodes of the graph replacing the gate are also connected to a pair of pins.

**Graph $G_i$, $2 \leq i \leq n$:** It has two calls to $G_{i-1}$ corresponding to the two calls to $C_{i-1}$ in $C_i$. For each of the explicit gates we replace it by a corresponding subgraph depending on whether it a AND or an OR gate. Again as in $G_1$ if the input of the gate is connected to pins then the corresponding in-nodes are connected to two pins.

An example of this construction appears in Figure **??**. The reader should notice that the construction produces a hierarchical description of the graph that would be obtained if the reduction of [KSS89] were applied on the circuit produced by the expansion $E(C)$ of the hierarchical specification $C$.

With the above observations, it is easy to see that the following lemmas from [KSS89] hold:

**Lemma 7.7** *The output of $C$ is 1 iff the graph $G$ has a 3-connected subgraph.*

**Lemma 7.8** *The above construction can be done in polynomial time.*

The theorem now follows from the above lemmas. ∎

The proofs of the following theorems also follow the same generic pattern as the proof of Theorem 7.6 above. The proof of Theorem 7.9 lifts the reduction in [KSS89] showing the P-hardness of approximating connectivity and the proof of Theorem 7.10 lifts the reduction in [AM86] showing the P-hardness of approximating the high degree subgraph problem.

**Theorem 7.9** *Unless $P = PSPACE$, the optimization version of the problem $\kappa$-HVC$_{HG}(G)$ and $\kappa$-HEC$_{HG}(G)$ cannot be approximated to within a factor of $c < 2$, even for simple strongly 1-level-restricted hierarchical specifications of $G$.*

**Theorem 7.10** *Unless $P = PSPACE$, the optimization version of the problem $HDSP_k$ cannot be approximated to within a factor $c < 2$ even for simple strongly 1-level-restricted hierarchical specifications of $G$.*

# 8 Conclusions and Related Work

We have presented polynomial time approximation algorithms with good performance guarantees for several natural PSPACE-complete problems for hierarchical specifications. We have also presented results concerning non-approximability of optimization version of the monotone circuit value problem, linear programming and high degree vertex and edge connectivity problems. Our proofs of non-approximability can be extended so as to apply to $O(\log \eta)$-bandwidth bounded hierarchical specifications, where $\eta$ is the size of the instance obtained after expanding the given specification. The question of whether the high degree subgraph and high connectivity problems for hierarchical specifications can be approximated to some constant factor of the optimal is open.

In [MRHR93] we have shown that efficient approximation algorithms can be obtained for hierarchically specified unit disk graphs. In [MHR93], we consider the complexity of finding polynomial time approximation schemes for hierarchically specified planar graphs. In [CF+93a, CF+93a] Condon et al. give a characterization of PSPACE in terms of probabilistically checkable debate systems and use this characterization to show that many natural PSPACE-hard problems cannot be approximated. Intriguingly enough, all the problems listed in Table 1 are known to have NC approximation algorithms when the problem instances are specified non-hierarchically [KW85, PSZ89]. Moreover, each of the problems shown to have a polynomial time optimal solution in [LW87a, Le88, Le89, Wi90] (eg. minimum spanning tree, planarity testing) when the problem is specified hierarchically, has an NC algorithm, when the problem instance is presented non-hierarchically. In [HM+93] we have shown that for every problem $\Pi$ in MAX SNP there is an NC approximation algorithm $A_\Pi$ with a constant performance guarantee. All the problems for which we have approximation algorithms in the hierarchical case belong to MAX SNP in the non-hierarchical case. While there are problems whose non-hierarchical versions can be solved in NC, but their hierarchical versions are PSPACE-hard [LW92], the results here and in [LW87a, Le88, Le89, Wi90] suggest that there is a strong relationship between a problem having an NC algorithm in the non-hierarchical case and a polynomial time algorithm in the hierarchical case. Understanding this relationship may well lead to a paradigm for translating known NC algorithms in the literature, for problems when specified non-hierarchically, to polynomial time algorithms for the same problems when the instances are specified hierarchically.

**Acknowledgements:** We thank Venkatesh Radhakrishnan and Richard Stearns for several constructive suggestions and Lefteris Kirousis for making available the journal version of [KSS89]. The first author also expresses his thanks to Egon Wanke for fruitful discussions on hierarchical specifications during his visit to Bonn.

**Algorithm FMAX 3SAT**

**Input:** A 3SAT formula $F$ and its associated bipartite graph.

1. Transform the bipartite graph $G$ corresponding to $F$ into a new bipartite graph $G'$ in which the we have one vertex for each variable, one vertex for each clause and if a clause $c_i = (x \vee y \vee z)$ then we have an edge from vertex corresponding to $c_i$ to the vertex corresponding to $x$.

   **Remark:** Step 1 intuitively breaks the original bipartite graph into stars with a variable node as the center of each star.

2. **For** each variable $x_i$, $1 \leq i \leq n$ **do**
   **Begin**

   (a) Compute the sets $PV^{x_i}$ and $NV^{x_i}$ defined as

   $$PV^{x_i} = \{w \mid w \text{ is a clause node } \mathbf{adjacent\ to}\ x_i \text{ in } G' \text{ and } x_i \text{ appears unnegated in } w\}$$

   $$NV^{x_i} = \{w \mid w \text{ is a clause node } \mathbf{adjacent\ to}\ x_i \text{ in } G' \text{ and } x_i \text{ appears negated in } w\}$$

   (b) **If** $|PV^{x_i}| \geq |NV^{x_i}|$ **then** set $x_i$ to true **else** set $x_i$ to false.

   **End**

3. **Output:** The satisfying assignment to the variables of $F$.

Figure 9: A Heuristic for Non-hierarchical Instances of MAX 3SAT

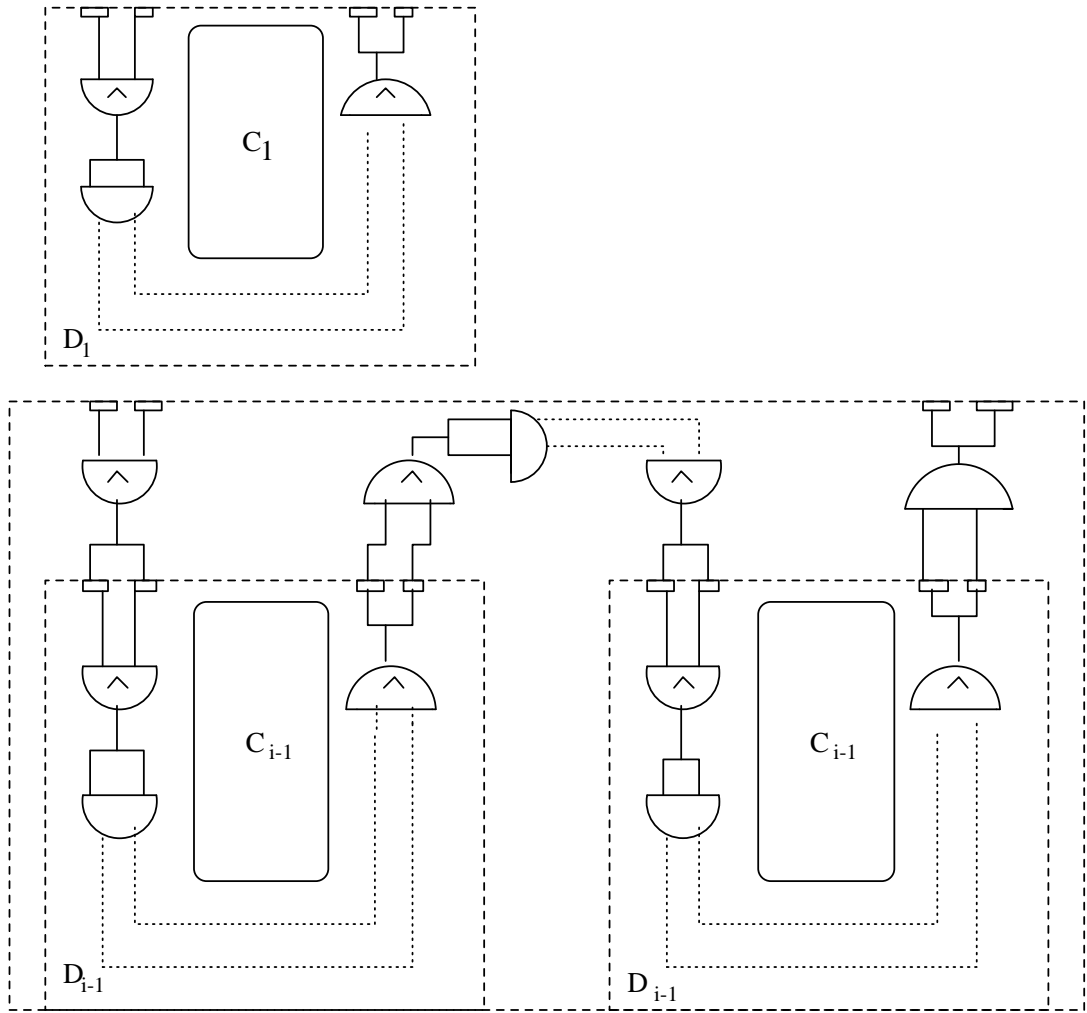

Figure 12: Construction of $D_i$, $1 \leq i < n$

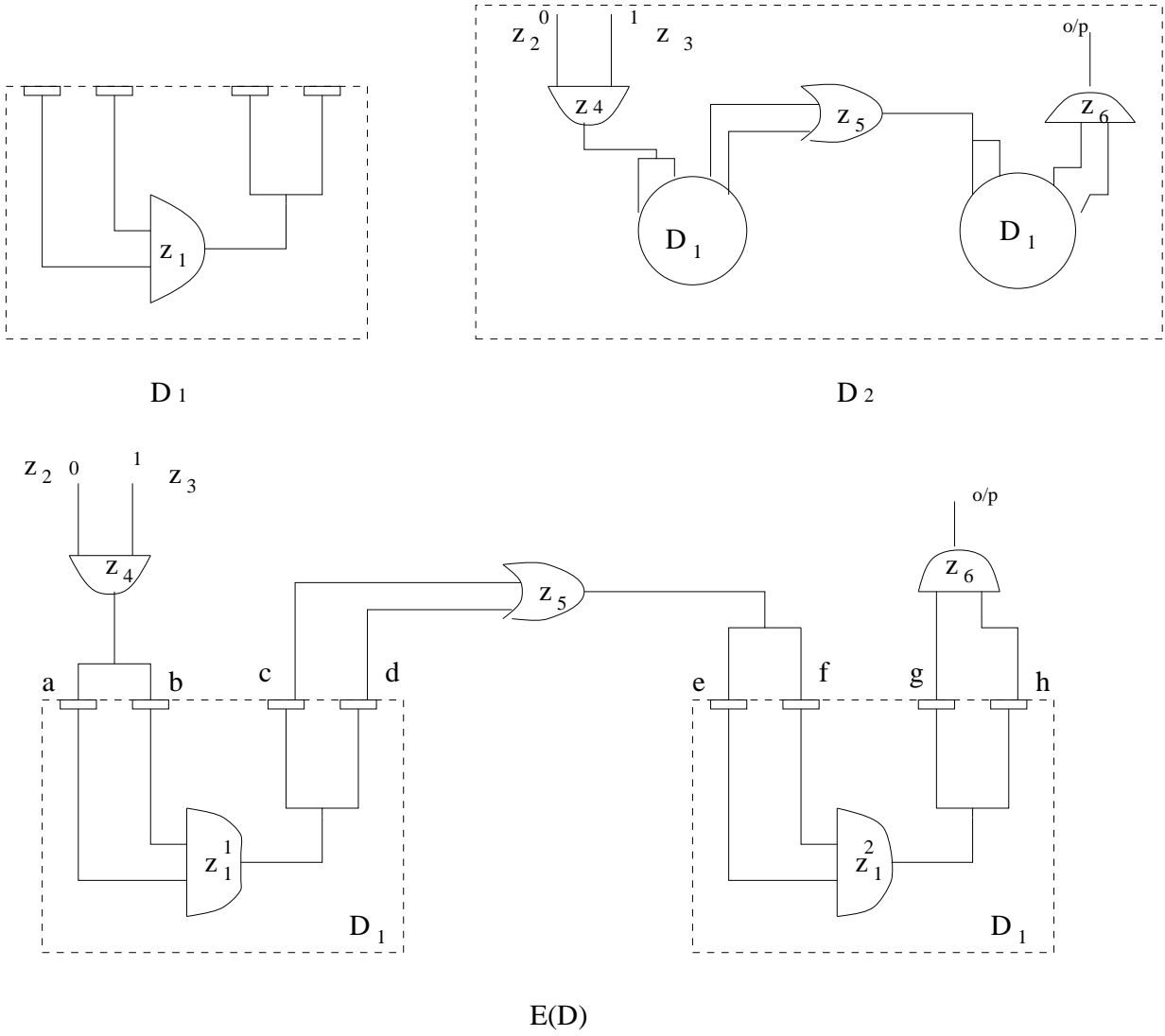

Figure 14: Example of a circuit represented hierarchically. $E(F)$ represents the actual circuit.